\documentclass[11pt, reqno]{amsart}

\usepackage[utf8]{inputenc}
\usepackage[T1]{fontenc}
\usepackage{microtype}
\usepackage{lmodern}
\usepackage{stmaryrd}
\usepackage{mathrsfs}

\usepackage{
    graphicx, 
    comment,
    fancyhdr, 
    amscd,
    amssymb,
    amsmath, 
    tikz, 
    tikz-cd, 
    mathtools,
    xfrac,
    booktabs,
    longtable,
    overpic,
    colortbl,
    lscape,
    xcolor
}

\usepackage[shortlabels]{enumitem}
\setlist[enumerate]{leftmargin=*}

\definecolor{vpurple}{HTML}{8E44AD}

\usepackage[
    colorlinks=true,
    linkcolor=blue,
    citecolor=blue,
    urlcolor=vpurple,
]{hyperref}

\usepackage[
    backend=biber,
    style=alphabetic,
    url=true,
]{biblatex}
\usetikzlibrary{calc}

\usetikzlibrary{decorations.pathmorphing}

\tikzset{curve/.style={settings={#1},to path={(\tikztostart)
    .. controls ($(\tikztostart)!\pv{pos}!(\tikztotarget)!\pv{height}!270:(\tikztotarget)$)
    and ($(\tikztostart)!1-\pv{pos}!(\tikztotarget)!\pv{height}!270:(\tikztotarget)$)
    .. (\tikztotarget)\tikztonodes}},
    settings/.code={\tikzset{quiver/.cd,#1}
        \def\pv##1{\pgfkeysvalueof{/tikz/quiver/##1}}},
    quiver/.cd,pos/.initial=0.35,height/.initial=0}

\tikzset{tail reversed/.code={\pgfsetarrowsstart{tikzcd to}}}
\tikzset{2tail/.code={\pgfsetarrowsstart{Implies[reversed]}}}
\tikzset{2tail reversed/.code={\pgfsetarrowsstart{Implies}}}

\tikzset{no body/.style={/tikz/dash pattern=on 0 off 1mm}}
\title{Combination Theorems and the Farrell--Jones Conjecture for Graph Products of Pairs}
\author{Ioannis Gkeneralis$^{\dagger}$}
\author{Stratos Prassidis$^{\ddagger}$}

\usepackage[margin=2.5cm,head=.25cm,foot=1.5cm]{geometry}

\numberwithin{equation}{section}

\let\oldtocsection=\tocsection
\let\oldtocsubsection=\tocsubsection
\let\oldtocsubsubsection=\tocsubsubsection
\renewcommand{\tocsection}[2]{\hspace{0em}\oldtocsection{#1}{#2}}
\renewcommand{\tocsubsection}[2]{\hspace{1em}\oldtocsubsection{#1}{#2}}
\renewcommand{\tocsubsubsection}[2]{\hspace{2em}\oldtocsubsubsection{#1}{#2}}

\newtheorem{theorem}{Theorem}[section]
\newtheorem{prop}[theorem]{Proposition}
\newtheorem{lemma}[theorem]{Lemma}
\newtheorem{corollary}[theorem]{Corollary}
\newtheorem{question}[theorem]{Question}

\theoremstyle{definition}

\newtheorem{definition}[theorem]{Definition}
\newtheorem{example}[theorem]{Example}
\newtheorem{remark}[theorem]{Remark}

\makeatletter
\newcommand{\proofstep}[1]{%
  \par% ensure starting on a new paragraph
  \addvspace{\smallskipamount}% some vertical space
  \textit{#1\@addpunct{.}}\enspace\ignorespaces
}
\makeatother

\definecolor{vblue}{HTML}{0089cf}
\definecolor{vorange}{HTML}{CF4600}
\definecolor{vgreen}{HTML}{1CC900}
\renewcommand{\colon}{
  \nobreak
  \mskip 2mu
  \mathpunct{}
  \nonscript
  \mkern-\thinmuskip
  \mathopen{:}
  \mskip 6mu plus 1mu\relax
} % fixes an issue where there is too much spacing between \colon and \operatorname

\makeatletter
\renewcommand{\paragraph}{%
  \@startsection{paragraph}{4}{\z@}%
  {1.75ex plus .5ex minus .2ex}% space before
  {.75ex}% space after; positive => heading on its own line
  {\normalfont\normalsize\bfseries}%
}
\makeatother

\begin{document}

\begin{abstract}
We study the Farrell--Jones conjecture for graph products of pairs and for fundamental groups of polyhedral products. We introduce a pair-deletion decomposition and a local notion of FJ-admissibility, and prove a vertexwise combination theorem for graph products of pairs. The criterion applies, in particular, to ordinary graph products, free-factor pairs, suitable normal pairs, and classes controlled by acylindrical Bass--Serre actions.

For polyhedral products, we combine these results with their description by relative graph products and graph products of pairs. We also obtain bounds and exact computations for amenable category, together with applications to simplicial volume and Borel rigidity of closed aspherical polyhedral-product manifolds.
\end{abstract}

\maketitle

${}^{\dagger}$Department of Mathematics, Aristotle University of Thessaloniki, Greece\\
Email: \href{mailto:igkeneralis@math.auth.gr}{\texttt{igkeneralis@math.auth.gr}}

\vspace{0.5em}

${}^{\ddagger}$Department of Mathematics, National Kapodistrian University of Athens, Greece\\
Email: \href{mailto:efprass@uoa.gr}{\texttt{efprass@uoa.gr}}

\tableofcontents

%--------------
\section{Introduction}

Polyhedral products form a flexible construction connecting the combinatorics of simplicial complexes with homotopy theory, geometric group theory, and manifold topology. Given a simplicial complex \(L\) and a collection of pairs of spaces \((\underline A,\underline B)=\{(A(i),B(i))\}_{i\in I}\), the associated polyhedral product \(Z_L(\underline A,\underline B)\) is assembled from products in which the coordinates indexed by a simplex of \(L\) are allowed to vary in the spaces \(A(i)\), while the remaining coordinates are constrained to lie in \(B(i)\). Classical examples include real, complex, and quaternionic moment-angle complexes, Davis--Januszkiewicz spaces, and classifying spaces of right-angled Artin and Coxeter groups; see, for instance, \cite{bbcg2009,bbcg2010,davis_asphericity}.

For aspherical polyhedral products, the interaction between the combinatorics of \(L\) and the fundamental groups of the defining pairs is particularly strong. Davis showed that the fundamental group can be described as a \emph{relative graph product}, determined by the actions \(G_i=\pi_1(A(i))\curvearrowright E_i=\pi_0\bigl(p_i^{-1}(B(i))\bigr)\), where \(p_i\colon\widetilde A(i)\to A(i)\) denotes the universal cover. When these actions are transitive, choosing a stabilizer \(H_i\leq G_i\) identifies \(E_i\cong G_i/H_i\), and the resulting fundamental group becomes the graph product of pairs \(\Pi_{L_1}(G_i,H_i)\) introduced by Januszkiewicz and Świ\k{a}tkowski \cite{JanuszkiewiczSwiatkowski}. Thus graph products of pairs occur naturally as fundamental groups of a broad class of polyhedral products and interpolate between ordinary graph products and the more general relative graph-product construction.

The purpose of this paper is to exploit this description in two directions. The first concerns invariants associated with families of subgroups, particularly amenable category and simplicial volume. The second concerns the Farrell--Jones conjecture and its rigidity consequences. A common feature of both directions is that global information about a graph product of pairs can be recovered from local information attached to the vertex pairs \((G_i,H_i)\).

\paragraph*{Main results}

The first result is a combination estimate for graph products of pairs. Let \(L\) be a finite flag complex and put \(G=\Pi_{L_1}(G_i,H_i)\). For every simplex \(J\in S(L)\), set
\[
G_J=
\prod_{j\in J}G_j
\times
\prod_{i\notin J}H_i,
\qquad
G_\varnothing=\prod_{i\in I}H_i.
\]

\begin{theorem}[Combination theorem]
\label{thm:intro-combination}
Let \(\mathcal F\) be a family of subgroups of \(G\). Then
\[
\operatorname{cat}_{\mathcal F}(G)
\leq
\max
\left\{
\operatorname{cat}_{\mathcal F|_{G_\varnothing}}
(G_\varnothing),
\;
\max_{\varnothing\neq J\in S(L)}
\left(
\operatorname{gd}(G_J)+|J|
\right)
\right\}.
\]
\end{theorem}

The same building-action argument applies more generally to relative graph products arising as fundamental groups of aspherical polyhedral products. In that setting, the relevant cell stabilizers are the residue stabilizers of the associated right-angled building; see Theorem~\ref{thm:polyhedral-product-combination} and the combination theorem of Capovilla--Li--L\"oh~\cite{capovilla_li_loh}.

The amenable case leads both to general upper bounds and, for split normal pairs, to exact computations.

\begin{theorem}[Amenable category]
\label{thm:intro-amenable}
Let \(L\) be a finite flag complex.
\begin{enumerate}[(i)]
\item
If \(G=\Pi_{L_1}(G_i,H_i)\) and every \(G_i\) is amenable, then \(\operatorname{cat}_{\mathrm{Am}}(G)\leq\dim L+1\).
\item
Suppose that \(G_i=H_i\rtimes Q_i\) for every \(i\), with each \(H_i\) amenable. Then \(\operatorname{cat}_{\mathrm{Am}}\Pi_{L_1}(G_i,H_i)=\operatorname{cat}_{\mathrm{Am}}\left(\mathop{\ast}_{L_1}Q_i\right)\).
\item
If, moreover, all \(Q_i=Q\), where \(Q\) is a non-trivial residually finite amenable group of type \(F\), then \(\operatorname{cat}_{\mathrm{Am}}\Pi_{L_1}(G_i,H_i)=\operatorname{vcd}(W_L)\).
\end{enumerate}
\end{theorem}

As a consequence, whenever \(M\) is an oriented closed connected \(m\)-dimensional manifold with \(\operatorname{cat}_{\mathrm{Am}}(\pi_1(M))<m\), Gromov's vanishing theorem gives \(\|M\|=0\). In particular, the resulting simplicial-volume applications do not require \(M\) itself to be aspherical.

Our second main result is a Farrell--Jones combination theorem for graph products of pairs.

\begin{theorem}[Farrell--Jones combination theorem]
\label{thm:intro-FJ}
Let \(\Gamma\) be a finite simplicial graph and let \(\{(G_v,H_v)\}_{v\in V(\Gamma)}\) be a family of group pairs. Assume that, for every \(v\in V(\Gamma)\), at least one of the following conditions holds:
\begin{enumerate}[(i)]
\item
\(H_v=\{1\}\) and \(G_v\in\mathcal C_{\mathrm{FJ}}\);
\item
\(H_v\trianglelefteq G_v\), \(G_v/H_v\in\mathcal C_{\mathrm{FJ}}\), and either \(G_v\in\mathcal C_{\mathrm{FJ}}\) or \(H_v\) is FJ-stable under cyclic extensions;
\item
\(G_v=H_v\ast K_v\) with \(H_v,K_v\in\mathcal C_{\mathrm{FJ}}\);
\item
\(G_v\in\mathcal C_{\mathrm{FJ}}\) and the pair \((G_v,H_v)\) satisfies the acylindricity criterion of Proposition~\ref{prop:acylindrical-FJ-admissible}.
\end{enumerate}
Then \(\Pi_\Gamma(G_v,H_v)\in\mathcal C_{\mathrm{FJ}}\).
\end{theorem}

The hypotheses are vertexwise: different vertices may be handled by different mechanisms. In particular, a single graph product of pairs may mix ordinary graph-product vertices, free-factor pairs, normal pairs controlled by cyclic-extension arguments, and non-normal pairs controlled by acylindrical Bass--Serre actions. Condition~{\rm (ii)} includes, for example, virtually solvable kernels and elementary amenable kernels of finite Hirsch length.

\begin{corollary}[Borel rigidity]
\label{cor:intro-Borel}
Let \(M=Z_L(\underline A,\underline B)\) be a closed aspherical manifold polyhedral product of dimension at least \(5\). Suppose that each local action \(G_i\curvearrowright E_i\) is transitive and that, after choosing \(H_i\leq G_i\) with \(E_i\cong G_i/H_i\), every pair \((G_i,H_i)\) satisfies one of the hypotheses of Theorem~\ref{thm:intro-FJ}. Then \(M\) satisfies the Borel conjecture.
\end{corollary}

The proofs of the two combination results use rather different mechanisms. For Theorem~\ref{thm:intro-combination}, the key input is the action of the fundamental group on the standard realization of the associated right-angled building. The stabilizers are residue groups, and in the transitive case they are precisely the groups \(G_J\) appearing above.

For Theorem~\ref{thm:intro-FJ}, we introduce a pair-deletion decomposition. Deleting a vertex \(v\) produces a retraction whose kernel is an amalgamated free product of copies of \(G_v\) over \(H_v\). The preimages of infinite cyclic subgroups are therefore of the form \(\left(\underset{s\in S}{\ast_{H_v}}G_v\right)\rtimes\mathbb Z\). This leads to the local notion of FJ-admissibility and reduces the global Farrell--Jones problem to a condition on the individual vertex pairs.

The two main directions are therefore complementary. The \(\mathcal F\)-category results exploit the geometry of the right-angled building and its cell stabilizers, while the Farrell--Jones results use the pair-deletion decomposition and inheritance through cyclic preimages. Both exhibit the same local-to-global principle: global properties of \(\Pi_{L_1}(G_i,H_i)\) are controlled by the local vertex pairs \((G_i,H_i)\).

The Farrell--Jones theorem also leads to topological rigidity applications. For closed aspherical manifold polyhedral products of dimension at least \(5\), the surgery-theoretic theorem of Bartels and L\"uck \cite{BartelsLuck} converts the Farrell--Jones conclusion into the Borel conjecture. We illustrate this with real moment-angle manifolds and M\"obius-type polyhedral products, and discuss surface-boundary pairs as a natural non-normal test case for the local theory.

\medskip

\paragraph*{Organization of the paper.}
In Section~\ref{sec:polyhedral_asphericity} we recall polyhedral products, the Davis--Kropholler asphericity criterion, and the associated right-angled building. Section~\ref{sec:fundamental-groups-relative-graph-products} discusses relative graph products, graph products of pairs, the fundamental group of a polyhedral product, and the corresponding building action and cell stabilizers. We then develop the combination estimates for families of subgroups and apply them to amenable category and simplicial volume (Sections~\ref{sec:combination_thms} and \ref{sec:applications}). Section~\ref{sec:classifying-spaces} specializes the preceding constructions to polyhedral products of classifying-space pairs \((BG_i,BH_i)\). Finally, in Section~\ref{sec:rigidity} we study the Farrell--Jones conjecture, introduce the pair-deletion construction and FJ-admissibility, prove the Farrell--Jones combination theorem for graph products of pairs, and derive the resulting Borel rigidity statements.

\paragraph*{Acknowledgements.}
The first author expresses his gratitude to P.~Capovilla for helpful discussions and suggestions on earlier versions of this paper.

\paragraph*{AI disclosure.} All the ideas, arguments, results and proofs were obtained without the use of AI. We used ChatGPT for proofreading. The authors take full responsibility for the results.

\section{Polyhedral products and asphericity}\label{sec:polyhedral_asphericity}

\subsection{Polyhedral products}

Let $L$ be a simplicial complex with finite vertex set \(I=\{1,\ldots,m\}\). We denote by $S(L)$ the face poset of $L$, including the empty simplex. For each $i\in I$, let $(A(i),B(i))$ be a pair of spaces with \(B(i)\subseteq A(i)\).

Throughout, we assume that each $A(i)$ is path connected and that each $B(i)$ is non-empty. When needed, we shall impose the additional hypotheses appearing in the asphericity results recalled below.

Given a simplex $J\in S(L)$, define
\[
Z_J(\underline A,\underline B)=\prod_{i\in J} A(i)\times\prod_{i\in I\setminus J} B(i),
\]
where
\[
(\underline A,\underline B) = \{(A(i),B(i))\}_{i\in I}.
\]
The \emph{polyhedral product} associated with $L$ and $(\underline A,\underline B)$ is the subspace
\[
Z_L(\underline A,\underline B) = \bigcup_{J\in S(L)}
Z_J(\underline A,\underline B)
\subseteq
\prod_{i\in I}A(i).
\]
If all pairs $(A(i),B(i))$ are equal to a fixed pair $(A,B)$, we simply write \(Z_L(A,B)\); for details see \cite{bbcg2009,bbcg2010}.

Equivalently, following Davis \cite{davis_asphericity}, for a point \(x=(x_i)_{i\in I}\in\prod_{i\in I}A(i)\) define its support by
\[
\operatorname{Supp}(x)=\{i\in I\mid x_i\in A(i)\setminus B(i)\}.
\]
Then
\[
Z_L(\underline A,\underline B) = \left\{ x\in\prod_{i\in I}A(i)
\ \middle|\
\operatorname{Supp}(x)\in S(L)
\right\}.
\]
Thus the combinatorics of $L$ determines which coordinates are allowed to leave the distinguished subspaces $B(i)$ simultaneously.

The same construction may also be viewed categorically. Regard $S(L)$ as a category under inclusion and consider the diagram
\[
J\longmapsto Z_J(\underline A,\underline B).
\]
Whenever $J\subseteq J'$, the inclusions $B(i)\subseteq A(i)$ induce a natural inclusion
\[
Z_J(\underline A,\underline B) \hookrightarrow Z_{J'}(\underline A,\underline B).
\]
The polyhedral product is therefore the colimit
\[
Z_L(\underline A,\underline B) \cong \operatorname*{colim}_{J\in S(L)} Z_J(\underline A,\underline B).
\]
This point of view, emphasized in the general theory of polyhedral products, makes explicit that the space is assembled functorially from the face poset of $L$ and the collection of pairs $(A(i),B(i))$.

We record several standard examples.

\begin{example}
If \((A(i),B(i))=([0,1],\{1\})\) for every $i\in I$, then \(K(L):=Z_L([0,1],\{1\})\) is the \emph{chamber associated with $L$} (for details see \cite[Ex. A.4.8, A.5.2]{davis_coxeter}). The order complex of the face poset $S(L)$ gives a canonical simplicial subdivision of $K(L)$. This chamber will appear later in the basic construction and in the description of the associated right-angled building.
\end{example}

\begin{example}
For the pair \((D^1,S^0)=([-1,1],\{-1,1\})\), the space \(\mathcal Z_L^{\mathbb R}=Z_L(D^1,S^0)\) is the \emph{real moment-angle complex}. When $L$ is the nerve of the facets of a simple polytope, this gives the usual real moment-angle manifold.
\end{example}

\begin{example}
For \((D^2,S^1)\), one obtains the classical moment-angle complex \(\mathcal Z_L=Z_L(D^2,S^1)\) (see \cite{tt}). Similarly, \(Z_L(D^4,S^3)\) is the quaternionic analogue of the moment-angle construction (see \cite{gken_rigidityCH_moment-angle}).
\end{example}

\begin{example}
If $G_i$ is a discrete group and $H_i\leq G_i$ is a subgroup, then the pairs \((BG_i,BH_i)\) give polyhedral products \(Z_L(\underline{BG},\underline{BH})\), whose fundamental groups will provide an important class of examples for the combination results developed later in the paper.
\end{example}

Two features of the definition will be important throughout. First, the construction depends functorially both on the simplicial complex and on the collection of pairs. Second, although the definition uses all faces of $L$, several fundamental-group constructions depend only on the $1$-skeleton of $L$. The interaction between the full simplicial structure and its $1$-skeleton becomes particularly important when $L$ is flag, and will be discussed in the next subsection.

\subsection{Flag complexes and right-angularity}

We now recall the combinatorial notions that underlie the asphericity theory of polyhedral products. Let $L$ be a simplicial complex with vertex set $I$, and denote by $L_1$ its $1$-skeleton.

\begin{definition}
A simplicial complex $L$ is called \emph{flag} if every finite subset $J\subseteq I$ whose vertices are pairwise joined by edges of $L$ spans a simplex of $L$.
Equivalently,
\[
J\in S(L)
\quad\Longleftrightarrow\quad
\{i,j\}\in S(L)
\text{ for every distinct }i,j\in J.
\]
\end{definition}

Thus a flag complex is completely determined by its $1$-skeleton. More precisely, given a simplicial graph $\Gamma$ with vertex set $I$, its \emph{flag completion}, or \emph{clique complex}, $\operatorname{Flag}(\Gamma)$ is the simplicial complex whose simplices are precisely the finite complete subgraphs of $\Gamma$. Hence
\[
L\text{ is flag}
\quad\Longleftrightarrow\quad
L=\operatorname{Flag}(L_1).
\]

Equivalently, $L$ is flag if every non-simplex contains a pair of vertices which does not span an edge. In terms of missing faces, this means that every minimal non-face of $L$ has cardinality two.

Throughout, we use the usual convention that a simplex $J$ with $|J|$ vertices has dimension \(\dim J=|J|-1\). Accordingly, a cube of type $J$ in the associated Davis realization has dimension $|J|$, and hence \(\dim\mathcal B_L=\dim L+1\).

\begin{example}
If $P$ is a poset, then its order complex \(\operatorname{Flag}(P)\) is the simplicial complex whose simplices are the non-empty finite chains \(p_0<p_1<\cdots<p_k\) in $P$. Since a finite collection of elements of $P$ which are pairwise comparable is totally ordered, $\operatorname{Flag}(P)$ is a flag complex. In particular, the barycentric subdivision of any simplicial or regular cell complex is flag.
\end{example}

The relevance of the flag condition to polyhedral products is already visible at the level of the $1$-skeleton. The fundamental groups that arise from polyhedral products are governed by commutation relations associated with the edges of $L$. If $L$ is flag, every collection of pairwise commuting directions corresponds to an actual simplex of $L$. Thus the simplicial structure contains no higher-dimensional \emph{missing faces} which are invisible from $L_1$. This principle will be central both in the asphericity criterion and in the description of the right-angled building associated with a polyhedral product.

We next recall the corresponding right-angled Coxeter group. Let $\Gamma=L_1$. The \emph{right-angled Coxeter group} associated with $\Gamma$ is
\[
W_\Gamma
=
\left\langle
s_i,\ i\in I
\ \middle|\
s_i^2=1,\;
[s_i,s_j]=1
\text{ whenever }\{i,j\}\in L_1
\right\rangle .
\]
When the underlying simplicial complex is clear, we shall simply write $W_L$.

For $J\subseteq I$, let \(W_J=\langle s_j\mid j\in J\rangle\) be the corresponding special subgroup. If $J$ spans a simplex of a flag complex $L$, then the generators in $J$ commute pairwise and \(W_J\cong (\mathbb Z/2)^{|J|}\).

Conversely, in a right-angled Coxeter group, $W_J$ is finite precisely when $J$ spans a complete subgraph of $L_1$. Hence, when $L$ is flag,
\[
J\in S(L)
\quad\Longleftrightarrow\quad
W_J\text{ is finite}.
\]
In the terminology of Coxeter groups, such subsets $J$ are called \emph{spherical}.

This correspondence identifies the face poset $S(L)$ of a flag complex with the poset of spherical subsets of the generating set of $W_L$. It is one of the basic reasons that flag complexes occur naturally in right-angled Coxeter geometry.

\subsection{Asphericity criteria of polyhedral products}
\label{subsec:asphericity}

We next recall the asphericity criterion for polyhedral products. The original criterion of M.Davis \cite{davis_asphericity} was subsequently refined by M.Davis and Kropholler \cite{davis_kropholler_corrigenda}; in particular, the corrected statement distinguishes those coordinates which contribute nontrivially to the associated right-angled structure. We formulate the corrected criterion below.

Throughout this subsection, let
\[
Z_L(\underline A,\underline B)=Z_L\bigl((A(i),B(i))_{i\in I}\bigr)
\]
be a polyhedral product, where the $A(i)$ and $B(i)$ are CW-complexes, each $A(i)$ is path connected, and a basepoint \(*_i\in B(i)\) has been fixed. For each $i\in I$, let \(p_i:\widetilde A(i)\longrightarrow A(i)\) be the universal covering and put \(\widetilde B(i)=p_i^{-1}(B(i))\). The group \(G_i=\pi_1(A(i),*_i)\) acts on the set of path components \(E_i:=\pi_0\bigl(\widetilde B(i)\bigr)\). Equivalently, $E_i$ may be identified with the relative $\pi_1$-set \(\pi_1(A(i),B(i),*_i)\). Following Davis--Kropholler \cite{davis_kropholler_corrigenda}, define \(I''=\{i\in I\mid |E_i|\neq 1\}\), and let \(L''=L_{I''}\) denote the full subcomplex of $L$ spanned by $I''$.

\begin{definition}
A vertex $i$ of a simplicial complex $K$ is called \emph{conelike} if it is joined by an edge to every other vertex of $K$.
\end{definition}

Thus, if $i$ is conelike and $K$ is flag, then \(K=\{i\}*K_{V(K)\setminus\{i\}}\). More generally, the collection of conelike vertices of a flag complex spans a simplex which splits off as a join factor. The following notion is needed in the corrected asphericity criterion.

\begin{definition}
Let $(A,B)$ be a pair with $A$ path connected, let \(\widetilde A\longrightarrow A\) be its universal covering, and let $\widetilde B$ be the inverse image of $B$.

The pair $(A,B)$ is called \emph{aspherical} if $\widetilde A$ is contractible and every path component of $\widetilde B$ is contractible.

The pair $(A,B)$ is called \emph{nearly aspherical} if $\widetilde A$ is contractible and every path component of $\widetilde B$ is acyclic.
\end{definition}

In particular, if each path component of $B$ is aspherical and the inclusion of every component induces an injection
\[
\pi_1(B,b)\longrightarrow \pi_1(A,b),
\]
then $(A,B)$ is an aspherical pair, and hence is nearly aspherical.

To state the general criterion, one further decomposition of the vertex set is required. Fix a field $\mathbb F$ and define disjoint subsets of $I\setminus I''$ by
\[
I_1=
\left\{
i\in I\setminus I''
\ \middle|\
\widetilde B(i)\text{ is acyclic over }\mathbb Z
\right\},
\]
and
\[
I_2=
\left\{
i\in I\setminus I''
\ \middle|\
\widetilde B(i)\text{ is not acyclic over }\mathbb F
\right\}.
\]
We shall write
\[
I'=I_2\cup I'',
\qquad
L'=L_{I'}.
\]

Following \cite{davis_kropholler_corrigenda}, we impose the auxiliary condition
\begin{equation}
\label{eq:DK-condition} I_1\sqcup I_2=I\setminus I''
\end{equation}
for some field $\mathbb F$. This excludes the exceptional situation in which the reduced homology of different $\widetilde B(i)$ is supported at incompatible primes; without such a hypothesis there is no equally simple necessary-and-sufficient formulation. We can now state the criterion.

\begin{theorem}[Davis--Kropholler]
\label{thm:Davis-Kropholler-asphericity}
Assume that \eqref{eq:DK-condition} holds. Then the polyhedral product \( Z_L(\underline A,\underline B)\) is aspherical if and only if the following conditions are satisfied:
\begin{enumerate}
    \item[\rm (i)]
    Each $A(i)$ is aspherical.

    \item[\rm (ii)]
    If $i\in I'$ is not conelike in $L'$, then $(A(i),B(i))$ is nearly aspherical.

    \item[\rm (iii)]
    Every vertex of $I_2$ is conelike in $L'$, the full subcomplex $\langle I_2\rangle$ is a simplex, and
    \[
    L'=\langle I_2\rangle *L''.
    \]

    \item[\rm (iv)]
    The simplicial complex $L'$ is flag.
\end{enumerate}
\end{theorem}

\begin{remark}
The distinction between $L$ and $L'$ is essential. In general, asphericity of the polyhedral product does \emph{not} require the whole complex $L$ to be flag. Coordinates for which $|E_i|=1$, or whose lifted distinguished subspace is homologically invisible in the sense above, may contribute cone-like directions which do not affect the homotopy type relevant to asphericity. The correct flag condition is therefore imposed on the reduced complex $L'$, rather than on $L$ itself.
\end{remark}

For the class of pairs most relevant to the sequel, the criterion simplifies considerably. Suppose that, for every $i\in I$ and every path component $C\subseteq B(i)$, the inclusion induces an injective homomorphism
\[
\pi_1(C)\longrightarrow \pi_1(A(i)),
\]
and suppose moreover that whenever $B(i)$ is path connected its image is a proper subgroup of $\pi_1(A(i))$. Then
\[
|E_i|\geq 2
\qquad\text{for every }i\in I,
\]
so that
\[
I''=I,\qquad L''=L, \qquad I_1=I_2=\varnothing.
\]
The corrected theorem therefore reduces to the following useful form.

\begin{corollary}
\label{cor:asphericity-standard}
Assume that for every $i\in I$:
\begin{enumerate}
    \item each path component of $B(i)$ is $\pi_1$-injective in $A(i)$;
    \item if $B(i)$ is path connected, then its fundamental group has proper image in $\pi_1(A(i))$.
\end{enumerate}
Then \(Z_L(\underline A,\underline B)\) is aspherical if and only if:
\begin{enumerate}
    \item each $A(i)$ is aspherical;
    \item for every nonconelike vertex $i$ of $L$, the pair $(A(i),B(i))$ is aspherical;
    \item $L$ is flag.
\end{enumerate}
\end{corollary}

\begin{proof}
Under the assumptions above, every $E_i$ has at least two  elements, and hence $I''=I$. Thus $L''=L$, while $I_1=I_2=\varnothing$, so Theorem \ref{thm:Davis-Kropholler-asphericity} reduces to conditions {\rm (i)}, {\rm (ii)}, and the flagness of $L$. Moreover, $\pi_1$-injectivity of the components of $B(i)$ identifies the components of $\widetilde B(i)$ with universal covers of the corresponding components of $B(i)$. Consequently, in this setting the aspherical-pair formulation gives the stated criterion.
\end{proof}

The role of the flag condition can be seen geometrically after passing to universal covers. Once the spaces $A(i)$ have been replaced by their universal covers, M.Davis \cite{davis_asphericity} reduces the relevant homotopy-theoretic problem to polyhedral products of pairs \((\operatorname{Cone}(E_i),E_i)\), where the $E_i$ are discrete sets of components of $\widetilde B(i)$. After discarding the indices for which $|E_i|=1$, the universal cover of this model is the standard realization of a right-angled building. When the corresponding simplicial complex is flag, this realization is contractible; indeed, it carries the standard $\operatorname{CAT}(0)$ cubical geometry.

This geometric model will be recalled in Subsection~\ref{subsec:basic-construction}. It will later provide the contractible $G$-CW-complex on which \(G=\pi_1 Z_L(\underline A,\underline B)\) acts, and whose cell stabilizers enter the combination theorem.

\subsection{The basic construction and right-angled buildings}
\label{subsec:basic-construction}

We next recall the basic construction associated with a mirrored space and explain its relation with right-angled Coxeter groups and right-angled buildings. This construction provides the geometric model that will later be used to study the action of the fundamental group of an aspherical polyhedral product.

Let $I$ be a finite set. An \emph{$I$-mirrored space} is a space $X$ together with a collection of subspaces \(\{X_i\}_{i\in I}\), called its \emph{mirrors}. For $x\in X$, define its type by \(I(x)=\{i\in I\mid x\in X_i\}\). For $J\subseteq I$, we put \(X_J=\bigcap_{j\in J}X_j\), with the convention $X_\varnothing=X$.

Let $(W,S)$ be a Coxeter system with its set of generators indexed by $I$, say \(S=\{s_i\}_{i\in I}\). For $J\subseteq I$, let \(W_J=\langle s_j\mid j\in J\rangle\) denote the corresponding special subgroup. The \emph{basic construction} associated with $(W,X)$ is the quotient \(\mathcal U(W,X)=(W\times X)/\sim\), where \((w,x)\sim(w',x')\) if and only if $x=x'$ and \(wW_{I(x)}=w'W_{I(x)}\). Equivalently,
\[
(w,x)\sim(w',x)
\quad\Longleftrightarrow\quad
w^{-1}w'\in W_{I(x)}.
\]
The group $W$ acts on $\mathcal U(W,X)$ by left multiplication on the first coordinate, \(g\cdot[w,x]=[gw,x]\). The image of $\{w\}\times X$ will be called the chamber $wX$. Thus $\mathcal U(W,X)$ is obtained by gluing together copies of the mirrored space $X$, with the gluing prescribed by the special subgroups $W_J$.

We apply this construction to the chamber determined by a simplicial complex. Let $L$ be a simplicial complex with vertex set $I$ and let $S(L)$ denote its face poset, including the empty simplex. Define \(K(L)=|\operatorname{Flag}(S(L))|\). The complex $K(L)$ is naturally an $I$-mirrored space. For each $i\in I$, the $i$-th mirror is
\[
K_i
=
\bigl|
\operatorname{Flag}
\{J\in S(L)\mid i\in J\}
\bigr|,
\]
and for a simplex $J\in S(L)$ one has
\[
K_J
=
\bigcap_{j\in J}K_j
=
\left|
\operatorname{Flag}(S(L)_{\geq J})
\right|,
\]
where
\[
S(L)_{\geq J}
=
\{J'\in S(L)\mid J\subseteq J'\}.
\]
Hence the intersections of mirrors are indexed by the simplices of $L$.

Let $W_L$ be the right-angled Coxeter group associated with the $1$-skeleton of $L$:
\[
W_L
=
\left\langle
s_i,\ i\in I
\ \middle|\
s_i^2=1,\;
[s_i,s_j]=1
\text{ whenever }\{i,j\}\in L_1
\right\rangle .
\]
The associated basic construction \(\Sigma_L :=\mathcal U(W_L,K(L))\) is the standard realization of the right-angled Coxeter system. When $L$ is flag, $\Sigma_L$ is the Davis complex of $W_L$.

The cells of $\Sigma_L$ may be described in terms of cosets of special subgroups. More precisely, simplices in the barycentric subdivision are indexed by chains \(J_0\subsetneq J_1\subsetneq\cdots\subsetneq J_k\) in the face poset $S(L)$, together with a coset of the corresponding special subgroup. In the cubical realization, a cube of type $J\in S(L)$ is represented by a coset \(wW_J\). Its stabilizer under the $W_L$-action is \(wW_Jw^{-1}\). Thus the face structure of $L$ determines both the cell types and the isotropy groups of the action.

When $L$ is flag, the cubical complex $\Sigma_L$ is $\operatorname{CAT}(0)$ and hence contractible. Equivalently, the link of every vertex is naturally isomorphic to $L$, and the flag condition is precisely Gromov's link condition for nonpositive curvature in this right-angled cubical setting. In particular,
\[
L\text{ flag}
\quad\Longrightarrow\quad
\mathcal U(W_L,K(L))
\text{ is contractible}.
\]

The Coxeter complex above should be viewed as an apartment in a more general right-angled building. We recall the form needed below. A \emph{chamber system over $I$} is a set $C$ together with equivalence relations \(\sim_i,\qquad i\in I\). Two chambers which are $\sim_i$-equivalent are said to be $i$-adjacent. For $J\subseteq I$, the equivalence relation generated by the $\sim_j$, $j\in J$, determines the $J$-residues.

Given such a chamber system $C$ and an $I$-mirrored space $X$, one defines analogously \(\mathcal U(C,X)=(C\times X)/\sim\), where \((c,x)\sim(c',x)\) whenever $c$ and $c'$ lie in the same $I(x)$-residue. The image of $\{c\}\times X$ is again called a chamber.

When $C$ is a right-angled building of type $(W_L,S)$ and $X=K(L)$, the space \(\mathcal U(C,K(L))\) is called the \emph{standard realization} of the building. Each apartment is naturally isomorphic to \(\mathcal U(W_L,K(L))\). Consequently, if $L$ is flag, the standard realization of the building is a $\operatorname{CAT}(0)$ cubical complex and in particular is contractible.

We now relate this construction to polyhedral products. Let \(E_i\) be a nonempty discrete set for each $i\in I$ and consider the family of pairs \((\operatorname{Cone}(E_i),E_i)\). The corresponding polyhedral product \(Z_L(\operatorname{Cone}(\underline E),\underline E)\) is a realization of the chamber system obtained from the product chamber system \(\prod_{i\in I}E_i\). After passing to the appropriate universal covering chamber system, Davis identifies its universal cover with a standard realization \(\mathcal U(C,K(L))\), where $C$ is a right-angled building of type $(W_L,S)$, provided the nontrivial factors $E_i$ are taken into account.

In particular, in the case
\(|E_i|>1\) for all \(i\in I\), and with $L$ flag, one obtains a contractible right-angled building realization (cf. \cite[Proposition~2.10 and Lemma~2.11]{davis_asphericity})
\[
\widetilde{ Z_L(\operatorname{Cone}(\underline E),\underline E)
} \cong \mathcal U(C,K(L)).
\]
This is the geometric model which enters Davis's proof of the asphericity theorem.

More generally, if some $E_i$ consists of a single point, then the corresponding coordinate is homotopically redundant. As recalled in Subsection~\ref{subsec:asphericity}, one passes to the full subcomplex
\[
L''=L_{I''}, \qquad I''=\{i\in I\mid |E_i|\neq1\}.
\]
The Davis--Kropholler correction \cite{davis_kropholler_corrigenda} shows that the relevant building is then the right-angled building associated with $L''$, and the appropriate contractibility condition is the flagness of $L''$ rather than necessarily of the whole complex $L$.

The importance of this construction for the sequel is that it converts the combinatorics of a simplicial complex into a contractible equivariant cell complex. For an aspherical polyhedral product \(Z_L(\underline A,\underline B)\), the fundamental group \(G=\pi_1 Z_L(\underline A,\underline B)\) will be described in the next section as a relative graph product. The Davis construction provides a natural contractible $G$-complex whose cell types are indexed by faces of $L$ and whose isotropy groups can be described in terms of the corresponding local subgroups. This is the geometric input needed later for the combination theorem of Capovilla--Li--Löh \cite{capovilla_li_loh}.

\section{Fundamental groups and relative graph products}\label{sec:fundamental-groups-relative-graph-products}

\subsection{Relative graph products}
\label{subsec:relative-graph-products}

We now recall the relative graph-product construction used by Davis to describe fundamental groups of polyhedral products. The ordinary graph product of groups is recovered as a special case, but the relative construction allows the local groups to act nontrivially on sets and is therefore adapted to the general pairs appearing in $Z_L(\underline A,\underline B)$.

Let $\Gamma$ be a simplicial graph with vertex set $I$. For each $i\in I$, let $G_i$ be a group acting on a nonempty set $E_i$. We write \(N_i=\ker(G_i\curvearrowright E_i)\) for the kernel of the action and \(\overline G_i=G_i/N_i\) for the corresponding faithful permutation group on $E_i$.

We first recall the ordinary graph product. The \emph{graph product} of the family $(G_i)_{i\in I}$ over $\Gamma$ is
\[
\mathop{\ast}_{\Gamma} G_i
=
\left(
\mathop{\ast}_{i\in I}G_i
\right)
\Big/
\left\langle\!\left\langle
[g_i,g_j]
\ \middle|\
g_i\in G_i,\ g_j\in G_j,\ \{i,j\}\in E(\Gamma)
\right\rangle\!\right\rangle .
\]
Thus elements belonging to vertex groups associated with adjacent vertices are required to commute.

For a complete subgraph $J\subseteq I$, the subgroup generated by the vertex groups $\{G_j\}_{j\in J}$ is naturally isomorphic to the direct product \(G_J=\prod_{j\in J}G_j\). When $L$ is flag and $\Gamma=L_1$, complete subgraphs of $\Gamma$ are precisely the simplices of $L$. Hence the simplices of $L$ index the standard product subgroups of the graph product.

The relative graph product is defined from the permutation data $(G_i,E_i)$. Put \(G=\prod_{i\in I}G_i\), \(N=\prod_{i\in I}N_i\), where \(N_i=\ker(G_i\curvearrowright E_i)\). The product group $G$ acts coordinatewise on \(Z_\Gamma(\operatorname{Cone}(\underline E),\underline E)\), and the kernel of this action is precisely $N$. Hence the effective group \(G/N\) acts on this polyhedral product.

Let \(\widetilde{ Z_\Gamma(\operatorname{Cone}(\underline E),\underline E) }\) be its universal cover, and let $\Gamma_\Gamma$ denote the group of all lifts of the $G/N$-action to this universal cover. Then there is a short exact sequence
\[
1\longrightarrow
\pi_1
Z_\Gamma(\operatorname{Cone}(\underline E),\underline E)
\longrightarrow
\Gamma_\Gamma
\longrightarrow
G/N
\longrightarrow1.
\]
Davis proves that $\Gamma_\Gamma$ depends only on the underlying graph.

The \emph{relative graph product} \(Z_\Gamma(\underline G,\underline E)\) is defined by the pullback diagram
\[
\begin{CD}
Z_\Gamma(\underline G,\underline E) @>>> \Gamma_\Gamma\\
@VVV @VVV\\
G @>>> G/N .
\end{CD}
\]
Consequently, it fits into an exact sequence
\[
1\longrightarrow N
\longrightarrow
Z_\Gamma(\underline G,\underline E)
\longrightarrow
\Gamma_\Gamma
\longrightarrow1.
\]

The construction interpolates between several familiar cases.

\begin{example}
If each action $G_i\curvearrowright E_i$ is faithful, then $N_i=1$ for all $i$, and the relative graph product is determined entirely by the induced permutation actions on the $E_i$.
\end{example}

\begin{example}
If every $E_i$ consists of a single point, then every local action is trivial and \(N_i=G_i\). In this degenerate case the right-angled chamber-system contribution disappears and the relative construction reduces to the corresponding product-type kernel.
\end{example}

\begin{example}
The ordinary graph product is recovered when the local permutation data are chosen so that no additional kernel correction is present. In particular, in the standard classifying-space case \((BG_i,*)\), the resulting fundamental group is the graph product of the groups $G_i$ over the $1$-skeleton of the indexing simplicial complex.
\end{example}

For a subset $J\subseteq I$, we shall write \(Z_J(\underline G,\underline E)\) for the relative graph product obtained by restricting the local data to the full subgraph on $J$. When $J$ spans a clique, this subgroup is the local group associated with the corresponding simplex. In the ordinary graph-product case it reduces to \(Z_J(\underline G,\underline E)\cong\prod_{j\in J}G_j\).

This local notation will be useful later when we describe stabilizers of cells in the right-angled building associated with an aspherical polyhedral product.

A particularly important case occurs when each action \(G_i\curvearrowright E_i\) is transitive. Choosing a point $e_i\in E_i$ and writing \(H_i=\operatorname{Stab}_{G_i}(e_i)\), we obtain an identification \(E_i\cong G_i/H_i\). In this case, M.Davis observes in \cite[Remark~2.15]{davis_asphericity} that the relative graph product \(Z_{L_1}(\underline G,\underline E)\) coincides with the graph product of pairs \(\Pi_{L_1}(G_i,H_i)\) introduced by Januszkiewicz and Świątkowski \cite{JanuszkiewiczSwiatkowski}. Thus graph products of pairs form a natural intermediate class between ordinary graph products and general relative graph products. In particular, the case $H_i=\{1\}$ for all $i$ recovers the ordinary graph product.

\begin{definition}\label{def:graph-prod-pairs}
Let \(\Gamma\) be a simplicial graph with vertex set \(I\), and let \(H_i\leq G_i\) for \(i\in I\). The \emph{graph product of pairs} \(\Pi_\Gamma(G_i,H_i)\) is the direct limit of the system \(G_J=\prod_{j\in J}G_j\times\prod_{i\notin J}H_i\), indexed by the cliques \(J\) of \(\Gamma\), with the natural inclusions.
\end{definition}

Finally, we emphasize that the relative graph product depends only on the graph $\Gamma$, and hence, in the polyhedral-product situation, only on the $1$-skeleton $L_1$. The higher-dimensional simplices of $L$ enter later through the flag condition and the geometry of the associated right-angled building.

\subsection{The fundamental group of a polyhedral product}
\label{subsec:fundamental-group-polyhedral-product}

We now relate the topology of a polyhedral product to the relative graph-product construction introduced in the previous subsection.

Let \((\underline A,\underline B)=\{(A(i),B(i))\}_{i\in I}\) be a family of pairs of CW-complexes indexed by the vertex set $I$ of a simplicial complex $L$. We assume throughout that every $A(i)$ is path connected. For each $i\in I$, let \(p_i:\widetilde A(i)\longrightarrow A(i)\) be the universal covering and put \(\widetilde B(i)=p_i^{-1}(B(i))\). Set \(G_i=\pi_1(A(i))\). The group $G_i$ acts by deck transformations on $\widetilde A(i)$ and hence permutes the path components of $\widetilde B(i)$. We denote the resulting $G_i$-set by \(E_i=\pi_0\bigl(\widetilde B(i)\bigr)\).

The collection \((\underline G,\underline E)=\{(G_i,E_i)\}_{i\in I}\) will be called the \emph{fundamental-group data} associated with $(\underline A,\underline B)$.

The first observation is that the fundamental group of a polyhedral product is insensitive to the simplices of dimension greater than one.

\begin{lemma}[M.Davis {cf.~\cite[Lemma 2.17]{davis_asphericity}}]
\label{lem:pi1-one-skeleton}
Let $L_1$ denote the $1$-skeleton of $L$. Then the inclusion \(Z_{L_1}(\underline A,\underline B) \longrightarrow Z_L(\underline A,\underline B)\) induces an isomorphism
\[
\pi_1\bigl(Z_{L_1}(\underline A,\underline B)\bigr)
\stackrel{\cong}{\longrightarrow}
\pi_1\bigl(Z_L(\underline A,\underline B)\bigr).
\]
In particular, \(\pi_1\bigl(Z_L(\underline A,\underline B)\bigr)\) depends only on the $1$-skeleton $L_1$.
\end{lemma}

\begin{remark}
This fact is conceptually important for what follows. The fundamental group is determined entirely by the graph $L_1$, whereas the asphericity of the polyhedral product depends, in general, on the full simplicial structure of $L$. In the flag case these two pieces of data are tied together by \(L=\operatorname{Flag}(L_1)\). Thus the same graph which governs the commutation relations in the fundamental group also determines the higher-dimensional simplices relevant to the right-angled building.
\end{remark}

The fundamental group can now be described explicitly in terms of the relative graph product of the local data.

\begin{theorem}[M.Davis {cf.~\cite[Theorem 2.18]{davis_asphericity}}]
\label{thm:Davis-fundamental-group}
Let $L$ be a simplicial complex with vertex set $I$, and let $(\underline A,\underline B)$ be a family of pairs as above. Let \((\underline G,\underline E) = \{(G_i,E_i)\}_{i\in I}\) be the associated fundamental-group data, where
\[
G_i=\pi_1(A(i)),
\qquad
E_i=\pi_0\bigl(\widetilde B(i)\bigr).
\]
Then there is a natural isomorphism \(\pi_1\bigl(Z_L(\underline A,\underline B)\bigr) \cong Z_{L_1}(\underline G,\underline E)\), where the right-hand side is the relative graph product over the $1$-skeleton $L_1$.
\end{theorem}

\begin{proof}
For each $i\in I$, collapse every path component of $\widetilde B(i)$ to a point. This gives a $G_i$-equivariant map \(\widetilde B(i)\longrightarrow E_i\). Since $\widetilde A(i)$ is path connected, this map extends to a $G_i$-equivariant map
\[
q_i:\widetilde A(i)\longrightarrow \operatorname{Cone}(E_i).
\]
Taking polyhedral products yields a map
\[
q: Z_L(\widetilde{\underline A},\widetilde{\underline B}) \longrightarrow Z_L(\operatorname{Cone}(\underline E),\underline E).
\]

The universal covering of \(Z_L(\widetilde{\underline A},\widetilde{\underline B})\) is obtained naturally by pulling back, along $q$, of the universal covering of \(Z_L(\operatorname{Cone}(\underline E),\underline E)\). Consequently, the deck-transformation group of the former is identified with the group of lifts determined by the latter. In detail, the pullback description therefore induces a homomorphism from the deck-transformation group of $\widetilde Z_L(\widetilde{\underline A},\widetilde{\underline B})$ to the group of lifts of the action on $\widetilde Z_L(\operatorname{Cone}(\underline E),\underline E)$. By the defining pullback property of the relative graph product, this homomorphism lifts to \(\pi_1Z_L(\underline A,\underline B) \longrightarrow Z_{L_1}(\underline G,\underline E)\), and the pullback construction shows that it is an isomorphism.

By the relative graph-product construction recalled in Subsection~\ref{subsec:relative-graph-products}, this group of lifts is precisely \(Z_{L_1}(\underline G,\underline E)\). Since the deck-transformation group of the universal covering of $Z_L(\underline A,\underline B)$ is its fundamental group, the desired isomorphism follows:
\[
\pi_1\bigl(Z_L(\underline A,\underline B)\bigr)
\cong
Z_{L_1}(\underline G,\underline E).
\]
\end{proof}

\begin{example}[Graph products of groups]
\label{ex:graph-products-BG}
Let \((A(i),B(i))=(BG_i,*)\). Then \(\widetilde A(i)=EG_i\) and the inverse image of the basepoint is a free $G_i$-orbit. Hence \(E_i\cong G_i\) as a $G_i$-set. Theorem~\ref{thm:Davis-fundamental-group} therefore recovers the ordinary graph product: \(\pi_1 Z_L(\underline{BG},*)\cong\mathop{\ast}_{L_1}G_i\).
\end{example}

\begin{example}[Right-angled Artin groups]
If \((A(i),B(i))=(S^1,*)\) for every $i\in I$, then $G_i\cong\mathbb Z$, and hence \(\pi_1 Z_L(S^1,*) \cong A_{L_1}\), the right-angled Artin group associated with the graph $L_1$.
\end{example}

\begin{example}[Right-angled Coxeter groups]
Let \((A(i),B(i))=(B\mathbb Z/2,*)\). Then we have that  \(\pi_1 Z_L(B\mathbb Z/2,*) \cong W_{L_1}\), the right-angled Coxeter group associated with $L_1$.
\end{example}

\begin{example}[Nontrivial relative data]
Suppose that each $B(i)$ has several lifts in $\widetilde A(i)$ and that the induced action \(G_i\curvearrowright E_i\) has nontrivial kernel \(N_i=\ker(G_i\curvearrowright E_i)\). Then the fundamental group is no longer described solely by the ordinary graph product of the $G_i$. Instead, \(\pi_1 Z_L(\underline A,\underline B) \cong Z_{L_1}(\underline G,\underline E)\), and the kernel data $N_i$ enter essentially into the relative graph-product construction.
\end{example}

\begin{remark}
The relative graph-product construction can be nontrivial even when all local groups $G_i$ are trivial. For example, if every $G_i$ is trivial and every $E_i$ is a two-point set, then \(Z_{L_1}(\underline G,\underline E)\cong\pi_1 Z_{L_1}(D^1,S^0)\). Thus the combinatorics of the chamber system itself can contribute nontrivial fundamental group.
\end{remark}

Theorem~\ref{thm:Davis-fundamental-group} identifies the acting group that will be relevant throughout the sequel. Under the asphericity hypotheses of Subsection~\ref{subsec:asphericity}, we shall write \(G=\pi_1 Z_L(\underline A,\underline B) \cong Z_{L_1}(\underline G,\underline E)\). The next step is to relate this relative graph-product description to the action of $G$ on the contractible right-angled building introduced in Subsection~\ref{subsec:basic-construction}, and in particular to determine the stabilizers of its cells.

\subsection{The building action and cell stabilizers}
\label{subsec:building-action}

We now combine the relative graph-product description of the fundamental group with the right-angled building introduced in Subsection~\ref{subsec:basic-construction}. The purpose of this subsection is to identify a natural contractible $G$-complex for \(G=\pi_1 Z_L(\underline A,\underline B)\) and to describe its isotropy groups in terms of residues of the associated building.

Throughout this subsection we first work in the clean situation in which
\[
|E_i|>1
\qquad\text{for every }i\in I,
\]
and $L$ is flag. The modifications required when some $E_i$ is a singleton are discussed at the end of the subsection.

Recall that \(G_i=\pi_1(A(i))\), \(E_i=\pi_0\bigl(\widetilde B(i)\bigr)\), and put \(\mathcal G=Z_{L_1}(\underline G,\underline E)\cong\pi_1 Z_L(\underline A,\underline B)\). For each $i\in I$, let \(N_i=\ker(G_i\curvearrowright E_i)\) and set \(N=\prod_{i\in I}N_i\).

Let \(C=\prod_{i\in I}E_i\) be the product chamber system. As recalled in Subsection~\ref{subsec:basic-construction}, there is a natural identification \(Z_L(\operatorname{Cone}(\underline E),\underline E)\cong\mathcal U(C,K(L))\). Let \(\widetilde{\mathcal B}_L = \widetilde{ Z_L(\operatorname{Cone}(\underline E),\underline E) }\). By \cite[Proposition~2.10]{davis_asphericity}, there exists a right-angled building $\widetilde C$ of type \((W_{L_1},S)\) such that \(\widetilde{\mathcal B}_L\cong\mathcal U(\widetilde C,K(L))\). Since $L$ is flag, this is the standard realization of $\widetilde C$ and is a $\operatorname{CAT}(0)$ cubical complex. In particular, \(\widetilde{\mathcal B}_L\) is contractible.

We next describe the action of $\mathcal G$ on this space.  Let $\Gamma_L$ denote the group of lifts of the effective $(\prod_iG_i)/N$-action on
\[
Z_L(\operatorname{Cone}(\underline E),\underline E)
\]
to its universal cover. By construction there is an exact sequence
\[
1
\longrightarrow
\pi_1
Z_L(\operatorname{Cone}(\underline E),\underline E)
\longrightarrow
\Gamma_L
\longrightarrow
\left(\prod_{i\in I}G_i\right)/N
\longrightarrow
1.
\]
The relative graph product $\mathcal G$ is obtained from this extension by the pullback construction recalled in Subsection~\ref{subsec:relative-graph-products}. In particular, there is a natural homomorphism \(\rho:\mathcal G\longrightarrow\Gamma_L, \) and hence an action \(\mathcal G\curvearrowright\widetilde{\mathcal B}_L\) given by \(g\cdot x=\rho(g)x\). Its kernel contains the subgroup $N$; in fact, with the above pullback description, \(\ker\rho=N\).

The action of $\Gamma_L$, and hence that of $\mathcal G$, preserves the type structure of the right-angled building. We therefore recall the relation between cell types and residues. For a chamber $c\in\widetilde C$ and a subset $J\subseteq I$, denote by \(\operatorname{Res}_J(c)\) the $J$-residue containing $c$.  If $J\in S(L)$, then this residue corresponds geometrically to the face of a chamber having type $J$. More precisely, consider a simplex of the barycentric chamber $K(L)$ corresponding to a chain \(J_0\subsetneq J_1\subsetneq\cdots\subsetneq J_k\) in $S(L)$. If $x$ lies in the relative interior of this simplex, then \(I(x)=J_0\). Consequently, in the basic construction \(\mathcal U(\widetilde C,K(L))\), the point $[c,x]$ depends only on the $J_0$-residue containing $c$.

The standard realization carries its natural cubical structure. Its cubes are indexed by residues of spherical type $J\in S(L)$. A cube of type $J$ has dimension $|J|$.

\begin{prop}
\label{prop:building-stabilizers}
Let $\sigma$ be a cube of type $J\in S(L)$ in $\widetilde{\mathcal B}_L$, and let $\operatorname{Res}_J(c)$ be the corresponding $J$-residue. Then \(\operatorname{Stab}_{\Gamma_L}(\sigma)=\operatorname{Stab}_{\Gamma_L}\bigl(\operatorname{Res}_J(c)\bigr)\). Consequently, \(\operatorname{Stab}_{\mathcal G}(\sigma)=\rho^{-1}\left(\operatorname{Stab}_{\Gamma_L}\bigl(\operatorname{Res}_J(c)\bigr) \right)\).
\end{prop}
\begin{proof}
Let $\sigma$ be a cube of type $J\in S(L)$, and let $\operatorname{Res}_J(c)$ be the corresponding $J$-residue.

By the definition of the standard realization \(\mathcal U(\widetilde C,K(L))\), a point in the relative interior of a cube of type $J$ depends only on the $J$-residue of the chamber containing it. Hence an element $\gamma\in\Gamma_L$ fixes $\sigma$ if and only if it preserves the corresponding $J$-residue setwise.

Indeed, if \(\gamma\cdot \sigma=\sigma\), then $\gamma$ sends any chamber containing $\sigma$ to another chamber containing the same cube, and therefore \(\gamma\operatorname{Res}_J(c)=\operatorname{Res}_J(c)\). Conversely, if \(\gamma\operatorname{Res}_J(c)=\operatorname{Res}_J(c)\), then $\gamma$ preserves the equivalence class defining the cube of type $J$, and hence fixes $\sigma$ setwise. Thus \(\operatorname{Stab}_{\Gamma_L}(\sigma)=\operatorname{Stab}_{\Gamma_L}\bigl(\operatorname{Res}_J(c)\bigr)\).

The action of $\mathcal G$ on $\widetilde{\mathcal B}_L$ factors through the homomorphism \(\rho:\mathcal G\longrightarrow\Gamma_L\). Therefore \(\operatorname{Stab}_{\mathcal G}(\sigma)=\rho^{-1}\left(\operatorname{Stab}_{\Gamma_L}(\sigma)\right)\), and using the first part we obtain \(\operatorname{Stab}_{\mathcal G}(\sigma)=\rho^{-1}\left(\operatorname{Stab}_{\Gamma_L}\bigl(\operatorname{Res}_J(c)\bigr)\right)\).

Finally, since \(N=\ker\rho\), every element of $N$ acts trivially on $\widetilde{\mathcal B}_L$. Hence \(N\subseteq \operatorname{Stab}_{\mathcal G}(\sigma)\) for every cube $\sigma$.
\end{proof}

\begin{remark}
Proposition~\ref{prop:building-stabilizers} is deliberately formulated in terms of residue stabilizers. In the general relative graph-product situation, the stabilizer depends not only on the simplex $J$, but also on the permutation data \(G_i\curvearrowright E_i\). Thus one should not identify it a priori with the direct product $\prod_{i\in J}G_i$.  The latter description is recovered in the ordinary graph-product case, as we record below.
\end{remark}

The preceding description simplifies considerably for ordinary graph products.

\begin{prop}
\label{prop:ordinary-graph-product-stabilizers}
Suppose that every $G_i$ acts simply transitively on $E_i$. Then \(E_i\cong G_i\) as a $G_i$-set, \(N_i=1\), and \(\mathcal G=\mathop{\ast}_{L_1}G_i\) is the ordinary graph product. For every simplex $J\in S(L)$, the stabilizer of a cell of type $J$ is conjugate to the standard parabolic subgroup \(G_J=\left\langle G_j\mid j\in J\right\rangle\). Since $L$ is flag, $J$ spans a complete subgraph of $L_1$, and hence \(G_J\cong\prod_{j\in J}G_j\). Therefore \(\operatorname{Stab}_{\mathcal G}(\sigma) \cong g\left(\prod_{j\in J}G_j \right) g^{-1}\) for some $g\in\mathcal G$.
\end{prop}

\begin{proof}
When the actions $G_i\curvearrowright E_i$ are simply transitive, Davis's relative graph product (cf. \cite[Remark~2.15]{davis_asphericity}) agrees with the ordinary graph product. The associated building is the standard right-angled building for this graph product. A cell of type $J$ is determined by a $J$-residue, and the stabilizer of the standard $J$-residue is the subgroup generated by the vertex groups $G_j$, $j\in J$.  Stabilizers of arbitrary $J$-residues are therefore conjugates of this subgroup. Since $J\in S(L)$ and $L$ is flag, its vertices are pairwise adjacent in $L_1$. The defining relations of the graph product consequently force the corresponding vertex groups to commute, giving \(G_J\cong\prod_{j\in J}G_j\).
\end{proof}

\begin{example}
\label{ex:BG-building-stabilizers}
For \((A(i),B(i))=(BG_i,*)\), the set $E_i$ is naturally a free and transitive $G_i$-set. Thus \(\pi_1Z_L(\underline{BG},*)\cong\mathop{\ast}_{L_1}G_i\), and the stabilizer of a cell of type $J$ in the associated right-angled building is conjugate to \(\prod_{j\in J}G_j\). This will provide the basic test case for the combination results in the next section.
\end{example}

We now give a more explicit description in the transitive case, using the graph products of pairs introduced by Januszkiewicz and Świątkowski in \cite{JanuszkiewiczSwiatkowski}.

Suppose that each $G_i$ acts transitively on $E_i$. Choose $e_i\in E_i$ and let \(H_i=\operatorname{Stab}_{G_i}(e_i)\). Then \(E_i\cong G_i/H_i\), and, by Davis, the relative graph product \(Z_{L_1}(\underline G,\underline E)\) coincides with the graph product of pairs \(\Pi_{L_1}(G_i,H_i)\) of Januszkiewicz--Świątkowski \cite{JanuszkiewiczSwiatkowski}. For $J\in S(L)$, define
\[
G_J
=
\prod_{j\in J}G_j
\times
\prod_{i\in I\setminus J}H_i.
\]
These groups form the local groups of the simple complex of groups associated with the graph product of pairs. Januszkiewicz and Świątkowski show that the corresponding development is the right-angled building associated with the graph product, and that the complex of groups arising from the group action on this building has precisely the local groups $G_J$.

\begin{prop}
\label{prop:pair-stabilizers}
Assume that every action
\(G_i\curvearrowright E_i\) is transitive, and let \(H_i=\operatorname{Stab}_{G_i}(e_i)\) for a chosen $e_i\in E_i$. Then \(G=Z_{L_1}(\underline G,\underline E)\cong\Pi_{L_1}(G_i,H_i)\). For every simplex $J\in S(L)$, the stabilizer of a cell of type $J$ in the associated right-angled building is conjugate to
\[
G_J=\prod_{j\in J}G_j\times\prod_{i\in I\setminus J}H_i.
\]
\end{prop}

\begin{proof}
Since the actions $G_i\curvearrowright E_i$ are transitive, we have \(E_i\cong G_i/H_i\). By
\cite[Remark~2.15]{davis_asphericity}, M.Davis identifies the corresponding relative graph product with the graph product of pairs \(\Pi_{L_1}(G_i,H_i)\) introduced by Januszkiewicz and Świątkowski \cite{JanuszkiewiczSwiatkowski}.

The graph product of pairs is the direct limit of the system of local groups
\[
G_J
=
\prod_{j\in J}G_j
\times
\prod_{i\in I\setminus J}H_i,
\qquad
J\in S(L).
\]
Januszkiewicz and Świątkowski construct the associated complex $D_G$ as the development of this simple complex of groups; see \cite[Proposition~2.1]{JanuszkiewiczSwiatkowski}. Moreover, the complex of groups determined by the $G$-action on $D_G$ is the original complex of groups with local groups $G_J$. Hence the stabilizer of the standard cell of type $J$ is $G_J$, and the stabilizer of any cell of type $J$ is conjugate to $G_J$, cf.~the proof of \cite[Theorem~3.1]{JanuszkiewiczSwiatkowski}.
\end{proof}

\begin{remark}
For $J=\varnothing$, Proposition~\ref{prop:pair-stabilizers} gives \(G_\varnothing=\prod_{i\in I}H_i\). Thus, unlike the ordinary graph-product case, vertex stabilizers need not be trivial.

On the other hand, if $H_i=\{1\}$ for every $i$, then \(G_J=\prod_{j\in J}G_j\), and Proposition~\ref{prop:pair-stabilizers} reduces to Proposition~\ref{prop:ordinary-graph-product-stabilizers}.
\end{remark}

\begin{remark}
Let \(N_i=\bigcap_{g\in G_i}gH_ig^{-1}\) be the core of $H_i$ in $G_i$. Then $N_i$ is precisely the kernel of the action \(G_i\curvearrowright G_i/H_i\). Consequently, \(N=\prod_{i\in I}N_i\) is contained in every local group $G_J$. This is consistent with the general description in Proposition~\ref{prop:building-stabilizers}, according to which the kernel of the action is contained in every cell stabilizer.
\end{remark}

For arbitrary, not necessarily transitive, $G_i$-sets $E_i$, the preceding graph-product-of-pairs description is no longer available in this form. In that general situation the simplex $J$ still determines a $J$-residue, but the associated stabilizer is naturally described by the residue-stabilizer formula of Proposition~\ref{prop:building-stabilizers}. Thus the relevant local groups for the sequel are \(\mathcal G_{c,J} :=\rho^{-1}\left( \operatorname{Stab}_{\Gamma_L} (\operatorname{Res}_J(c))\right)\). With this notation, \(\operatorname{Stab}_{\mathcal G}(\sigma)=\mathcal G_{c,J}\) for every cell $\sigma$ of type $J$ based at the chamber $c$.

Although $\mathcal G_{c,J}$ may depend on the orbit of the residue, its conjugacy class is constant on every $\mathcal G$-orbit of cells. Consequently, after choosing one representative from each $\mathcal G$-orbit of cells, the collection of groups \(\{\mathcal G_{c,J}\}\) is precisely the isotropy data required by combination theorems for groups acting on contractible CW-complexes.

Finally, suppose that some of the sets $E_i$ consist of a single point. As in Subsection~\ref{subsec:asphericity}, put
\[
I''=\{i\in I\mid |E_i|>1\},
\qquad
L''=L_{I''}. \]
The nontrivial right-angled building is then associated with $L''$. The preceding discussion applies after replacing $L$ by $L''$ and restricting the chamber data to the indices in $I''$; the remaining coordinates contribute to the kernel and product factors described by the Davis--Kropholler reduction.

We have therefore obtained a contractible $\mathcal G$-CW-complex \(\widetilde{\mathcal B}_L\) together with an explicit description of its cell stabilizers in terms of residue stabilizers. This is exactly the form of equivariant data to which the combination results recalled in the next section apply.

\section{Combination theorems for polyhedral-product groups}\label{sec:combination_thms}

\subsection{Families of subgroups and $\mathcal F$-category}
\label{subsec:F-category}

We begin by recalling the notion of $\mathcal F$-category introduced in the setting of open covers with fundamental-group constraints.

Let $G$ be a group. A \emph{family of subgroups} of $G$ is a nonempty collection $\mathcal F$ of subgroups of $G$ which is closed under conjugation and under passage to subgroups. Standard examples include
\[
\mathrm{Tr},\qquad
\mathrm{Fin},\qquad
\mathrm{Am},
\]
the families consisting, respectively, of the trivial subgroup, the finite subgroups, and the amenable subgroups of $G$.

Let $Y$ be a path-connected pointed space and let \(G=\pi_1(Y)\). Following Capovilla--Li--L\"oh \cite[Definition 1.2]{capovilla_li_loh}, the \emph{$\mathcal F$-category} of $Y$, denoted \(\operatorname{cat}_{\mathcal F}(Y)\), is the infimum of all $n\in\mathbb N$ for which there exists an open cover \(Y=U_0\cup\cdots\cup U_n\) such that, for every $i\in\{0,\ldots,n\}$ and every point $y\in U_i$, the image of the homomorphism \(\pi_1(U_i,y)\longrightarrow\pi_1(Y,y)\) belongs to $\mathcal F$.

For a group $G$, one defines \(\operatorname{cat}_{\mathcal F}(G) := \operatorname{cat}_{\mathcal F}(BG)\). This is independent of the choice of the classifying space $BG$. Moreover, if \(Y\longrightarrow Z\) induces an isomorphism on fundamental groups, then \(\operatorname{cat}_{\mathcal F}(Y)\leq\operatorname{cat}_{\mathcal F}(Z)\). In particular, \(\operatorname{cat}_{\mathcal F}(Y)\leq\operatorname{cat}_{\mathcal F}(\pi_1(Y))\).

If \(\mathcal E\subseteq\mathcal F\) are families of subgroups of $G$, then clearly \(\operatorname{cat}_{\mathcal F}(Y)\leq
\operatorname{cat}_{\mathcal E}(Y)\). Also,
\[
\operatorname{cat}_{\mathcal F}(G)=0
\quad\Longleftrightarrow\quad
G\in\mathcal F.
\]

For later use, if $H\leq G$, we denote by \(\mathcal F|_H=\{F\leq H\mid F\in\mathcal F\} \) the family obtained by restricting $\mathcal F$ to $H$.

The $\mathcal F$-category admits an equivalent equivariant formulation which is particularly well suited to the study of groups acting on contractible CW-complexes. Let $X$ be a $G$-CW-complex. We define \(\operatorname{cat}^{G}_{\mathcal F}(X)\) to be the infimum of all $n\in\mathbb N$ for which there exist an $n$-dimensional $G$-CW-complex $N$ with isotropy in $\mathcal F$ and a $G$-map \(X\longrightarrow N\). Here \emph{with isotropy in $\mathcal F$} means that the stabilizer of every cell of $N$ belongs to $\mathcal F$.

The equivariant invariant is monotone under $G$-maps: if \(X_1\longrightarrow X_2\) is a $G$-map, then \(\operatorname{cat}^{G}_{\mathcal F}(X_1)\leq\operatorname{cat}^{G}_{\mathcal F}(X_2)\). In particular, \(\operatorname{cat}^{G}_{\mathcal F}\) is invariant under $G$-homotopy equivalence.

The connection with the nonequivariant definition is the following.

\begin{prop}[Capovilla--Li--L\"oh {cf.~\cite{capovilla_li_loh}}]
\label{prop:F-category-universal-cover}
Let $Y$ be a connected CW-complex with fundamental group \(G=\pi_1(Y)\), and let \(\widetilde Y\) be its universal covering, endowed with the induced $G$-CW-structure. Then \(\operatorname{cat}_{\mathcal F}(Y)=\operatorname{cat}^{G}_{\mathcal F}(\widetilde Y)\).
In particular, \(\operatorname{cat}_{\mathcal F}(G)=\operatorname{cat}^{G}_{\mathcal F}(EG)\).
\end{prop}

This reformulation is the key point for our purposes. Indeed, if \(X\) is a contractible $G$-CW-complex, then $X$ may be used as a model on which to estimate the $\mathcal F$-category of $G$ by means of its cell stabilizers. In the polyhedral-product setting developed above, the right-angled building \(\widetilde{\mathcal B}_L\) provides precisely such a contractible $G$-complex.

We record several special cases. For the trivial family \(\mathrm{Tr}=\{\{1\}\}\), one has \(\operatorname{cat}_{\mathrm{Tr}}(G)=\operatorname{cd}(G)\). For the family $\mathrm{Fin}$ of finite subgroups, the corresponding invariant is closely related to virtual cohomological dimension. For the family $\mathrm{Am}$ of amenable subgroups, \(\operatorname{cat}_{\mathrm{Am}}(G)\) controls several homological invariants, including bounded cohomology and simplicial-volume-type invariants.

The usefulness of the equivariant formulation is that it permits one to estimate global invariants of $G$ from the isotropy groups of a contractible $G$-CW-complex. The required combination theorem will be recalled in the next subsection and then applied to the building action of \(G=\pi_1 Z_L(\underline A,\underline B)\).

\subsection{Combination via contractible $G$-CW-complexes}
\label{subsec:combination-G-CW}

We now recall the combination theorem of Capovilla--Li--L\"oh \cite{capovilla_li_loh}, which provides the group-theoretic input for our applications to polyhedral products.

Let $G$ be a group, let $\mathcal F$ be a family of subgroups of $G$, and let $X$ be a $G$-CW-complex. For a cell $\sigma\subseteq X$, we write \(G_\sigma=\operatorname{Stab}_G(\sigma)\). For each $i\geq 0$, let $\Sigma_i$ be a set of representatives of the $G$-orbits of $i$-cells of $X$. 

The basic principle is that, when $X$ is contractible, the $\mathcal F$-category of $G$ can be estimated from the corresponding invariants of the cell stabilizers. There are two complementary ways of doing this. We write \(\operatorname{gd}(G)\) for the minimal dimension of a model for \(BG\) (equivalently, of a free contractible \(G\)-CW-complex). One uses the geometric dimensions \(\operatorname{gd}(G_\sigma)\), whereas the other uses the restricted categories \(\operatorname{cat}_{\mathcal F|_{G_\sigma}}(G_\sigma)\). The general combination theorem permits the two estimates to be mixed in different dimensions.

\begin{theorem}[Capovilla--Li--L\"oh]
\label{thm:CLL-combination}
Let $G$ be a group, let $\mathcal F$ be a family of subgroups of $G$, and let $X$ be an $n$-dimensional $G$-CW-complex. For every $i\in\{0,\ldots,n\}$, let $\Sigma_i$ be a set of representatives of the $G$-orbits of $i$-cells of $X$.

Fix a subset \(\mathcal I\subseteq\{1,\ldots,n\}\). Define \(d_0 = \sup_{v\in\Sigma_0}\operatorname{cat}_{\mathcal F|_{G_v}}(G_v)\), and, inductively for $1\leq i\leq n$, set
\[
d_i
=
\begin{cases}
\displaystyle
\max\left\{
d_{i-1},
\sup_{\sigma\in\Sigma_i}
\bigl(
\operatorname{gd}(G_\sigma)+i
\bigr)
\right\},
&
i\in\mathcal I,
\\[4mm]
\displaystyle
d_{i-1}
+
\sup_{\sigma\in\Sigma_i}
\left(
\operatorname{cat}_{\mathcal F|_{G_\sigma}}(G_\sigma)+1
\right),
&
i\notin\mathcal I.
\end{cases}
\]
Then \(\operatorname{cat}_{\mathcal F}^{G}(EG\times X) \leq d_n\), where $G$ acts diagonally on $EG\times X$.

If, in addition, $X$ is non-equivariantly contractible, then \(\operatorname{cat}_{\mathcal F}(G)
\leq d_n\).
\end{theorem}

The strength of Theorem~\ref{thm:CLL-combination} lies in the freedom to choose $\mathcal I$. In each positive dimension one may use either geometric dimension or $\mathcal F$-category, depending on which produces the sharper estimate for the corresponding isotropy groups.

We record the two extreme choices separately.

\begin{corollary}
\label{cor:CLL-gd-bound}
Let $G$ act cellularly on a finite-dimensional, non-equivariantly contractible $G$-CW-complex $X$. Let $\Sigma_0$ be a set of $G$-orbit representatives of the vertices of $X$, and let $\Sigma_{\geq 1}$ be a set of $G$-orbit representatives of all positive-dimensional cells. Then
\[
\operatorname{cat}_{\mathcal F}(G)
\leq
\max
\left\{
\sup_{v\in\Sigma_0}
\operatorname{cat}_{\mathcal F|_{G_v}}(G_v),
\;
\sup_{\sigma\in\Sigma_{\geq1}}
\bigl(
\operatorname{gd}(G_\sigma)+\dim\sigma
\bigr)
\right\}.
\]
\end{corollary}

\begin{proof}
Apply Theorem~\ref{thm:CLL-combination} with \(\mathcal I=\{1,\ldots,\dim X\}\).
\end{proof}

At the opposite extreme one obtains an estimate involving only $\mathcal F$-categories of stabilizers.

\begin{corollary}
\label{cor:CLL-F-bound}
Let $G$ act cellularly on an $n$-dimensional, non-equivariantly contractible $G$-CW-complex $X$. Then
\[
\operatorname{cat}_{\mathcal F}(G)
\leq
\sup_{v\in\Sigma_0}
\operatorname{cat}_{\mathcal F|_{G_v}}(G_v)
+
\sum_{i=1}^{n}
\sup_{\sigma\in\Sigma_i}
\left(
\operatorname{cat}_{\mathcal F|_{G_\sigma}}(G_\sigma)+1
\right).
\]
\end{corollary}

\begin{proof}
Apply Theorem~\ref{thm:CLL-combination} with \(\mathcal I=\varnothing\).
\end{proof}

\begin{remark}
Neither of the preceding estimates dominates the other in general. The geometric-dimension bound is particularly effective when the positive-dimensional cell stabilizers have small geometric dimension, whereas the second estimate may be sharper when their $\mathcal F$-categories are substantially smaller than their geometric dimensions. The mixed form of Theorem~\ref{thm:CLL-combination} allows these two phenomena to be combined.
\end{remark}

The relevance of these results to polyhedral products is now immediate at the level of principle. In Subsection~\ref{subsec:building-action}, we associated to an aspherical polyhedral product a contractible right-angled building equipped with an action of its fundamental group and described the cell stabilizers in terms of residue stabilizers. We now insert this geometric model into Theorem~\ref{thm:CLL-combination}.

\subsection{The main combination theorem}
\label{subsec:main-combination}

We now apply the preceding combination principle to fundamental groups of aspherical polyhedral products.

Let \(Z=Z_L(\underline A,\underline B)\) be a polyhedral product satisfying the asphericity hypotheses of Subsection~\ref{subsec:asphericity}. For the moment, we work under the additional simplifying assumptions that \(|E_i|>1\) for all \(i\in I\), and that $L$ is flag. Put \(G=\pi_1(Z) \cong Z_{L_1}(\underline G,\underline E)\), where \(G_i=\pi_1(A(i))\), and \(E_i=\pi_0(\widetilde B(i))\).

By Subsection~\ref{subsec:building-action}, $G$ acts on the standard realization \(\mathcal B_L=\mathcal U(\widetilde C,K(L))\) of the associated right-angled building. Since $L$ is flag, $\mathcal B_L$ is a finite-dimensional $\operatorname{CAT}(0)$ cubical complex and hence is contractible.

Let \(n=\dim\mathcal B_L=\dim L+1\). For each $i\in\{0,\ldots,n\}$, choose a set $\Sigma_i$ of representatives of the $G$-orbits of $i$-dimensional cubes of $\mathcal B_L$.

If $\sigma\in\Sigma_i$, then $\sigma$ has a type \(J(\sigma)\in S(L)\), \(|J(\sigma)|=i\). Choosing a chamber $c_\sigma$ containing $\sigma$, put \(H_\sigma:=\mathcal G_{c_\sigma,J(\sigma)}=\operatorname{Stab}_G(\sigma)\). By Proposition~\ref{prop:building-stabilizers}, this group may equivalently be described as \(H_\sigma=\rho^{-1}\left(\operatorname{Stab}_{\Gamma_L}\bigl(\operatorname{Res}_{J(\sigma)}(c_\sigma) \bigr) \right)\).

We can now state the main combination theorem.

\begin{theorem}
\label{thm:polyhedral-product-combination}
Let \(Z=Z_L(\underline A,\underline B)\) be an aspherical polyhedral product as above, and let \(G=\pi_1(Z)\). Let $\mathcal F$ be a family of subgroups of $G$, and let \(\mathcal I\subseteq\{1,\ldots,n\}\), \(n=\dim L+1\). Define
\[
d_0
=
\sup_{\sigma\in\Sigma_0}
\operatorname{cat}_{\mathcal F|_{H_\sigma}}(H_\sigma),
\]
and, for $1\leq i\leq n$, define inductively
\[
d_i
=
\begin{cases}
\displaystyle
\max\left\{
d_{i-1},
\sup_{\sigma\in\Sigma_i}
\bigl(
\operatorname{gd}(H_\sigma)+i
\bigr)
\right\},
&
i\in\mathcal I,
\\[4mm]
\displaystyle
d_{i-1}
+
\sup_{\sigma\in\Sigma_i}
\left(
\operatorname{cat}_{\mathcal F|_{H_\sigma}}(H_\sigma)+1
\right),
&
i\notin\mathcal I.
\end{cases}
\]
Then \(\operatorname{cat}_{\mathcal F}\left(\pi_1 Z_L(\underline A,\underline B)\right)\leq d_n\).
\end{theorem}

\begin{proof}
By Subsection~\ref{subsec:building-action}, the group \(G=\pi_1Z_L(\underline A,\underline B)\) acts cellularly on the right-angled building realization \(\mathcal B_L=\mathcal U(\widetilde C,K(L))\). Since $L$ is flag, $\mathcal B_L$ is $\operatorname{CAT}(0)$ and therefore contractible.

For every cell $\sigma\subseteq\mathcal B_L$, its stabilizer is precisely the residue group \(G_\sigma=H_\sigma=\mathcal G_{c_\sigma,J(\sigma)}\) described in Proposition~\ref{prop:building-stabilizers}. Theorem~\ref{thm:CLL-combination} applied to the contractible $G$-CW-complex $\mathcal B_L$ therefore gives \(\operatorname{cat}_{\mathcal F}(G)\leq d_n\). Using \(G\cong\pi_1Z_L(\underline A,\underline B)\) gives the result.
\end{proof}

The two extreme choices of $\mathcal I$ give useful forms of Theorem~\ref{thm:polyhedral-product-combination}.

\begin{corollary}
\label{cor:polyhedral-product-gd-combination}
Under the assumptions of Theorem~\ref{thm:polyhedral-product-combination}, one has
\[
\operatorname{cat}_{\mathcal F}(G)
\leq
\max
\left\{
\sup_{\sigma\in\Sigma_0}
\operatorname{cat}_{\mathcal F|_{H_\sigma}}(H_\sigma),
\;
\sup_{\substack{\sigma\subseteq\mathcal B_L\\
                 \dim\sigma>0}}
\left(
\operatorname{gd}(H_\sigma)+\dim\sigma
\right)
\right\}.
\]
\end{corollary}

\begin{proof}
Take \(\mathcal I=\{1,\ldots,n\}\) in Theorem~\ref{thm:polyhedral-product-combination}.
\end{proof}

\begin{corollary}
\label{cor:polyhedral-product-F-combination}
Under the same assumptions,
\[
\operatorname{cat}_{\mathcal F}(G)
\leq
\sup_{\sigma\in\Sigma_0}
\operatorname{cat}_{\mathcal F|_{H_\sigma}}(H_\sigma)
+
\sum_{i=1}^{n}
\sup_{\sigma\in\Sigma_i}
\left(
\operatorname{cat}_{\mathcal F|_{H_\sigma}}(H_\sigma)+1
\right).
\]
\end{corollary}

\begin{proof}
Take $\mathcal I=\varnothing$ in Theorem~\ref{thm:polyhedral-product-combination}.
\end{proof}

\subsection{Graph products of pairs}
\label{subsec:combination-pairs}

We now specialize the preceding combination theorem to the transitive relative case. Besides producing explicit bounds, this connects the polyhedral-product construction with the graph products of pairs of Januszkiewicz and Świątkowski \cite{JanuszkiewiczSwiatkowski}.

Assume that \(E_i\cong G_i/H_i\) for every $i\in I$, and put \(G=\Pi_{L_1}(G_i,H_i)\cong Z_{L_1}(\underline G,\underline E)\). For every simplex $J\in S(L)$, define \(G_J=\prod_{j\in J}G_j\times\prod_{i\in I\setminus J}H_i\). By Proposition~\ref{prop:pair-stabilizers}, the stabilizer of a cube of type $J$ in the associated right-angled building is conjugate to $G_J$.

We therefore obtain the following explicit form of the combination estimate.

\begin{theorem}
\label{thm:pair-combination}
Let $L$ be a finite flag complex and let \(G=\Pi_{L_1}(G_i,H_i)\) be a graph product of pairs. Let $\mathcal F$ be a family of subgroups of $G$. Then
\[
\operatorname{cat}_{\mathcal F}(G)
\leq
\max
\left\{
\operatorname{cat}_{\mathcal F|_{G_\varnothing}}
(G_\varnothing),
\;
\max_{\varnothing\neq J\in S(L)}
\left(
\operatorname{gd}(G_J)+|J|
\right)
\right\},
\]
where \(G_\varnothing=\prod_{i\in I}H_i\) and \(G_J=\prod_{j\in J}G_j \times\prod_{i\notin J}H_i\). 
\end{theorem}

\begin{proof}
Let $D_G$ denote the standard right-angled building associated with the graph product of pairs. Since $L$ is flag, its standard cubical realization is contractible.

By Proposition~\ref{prop:pair-stabilizers}, a cube of type $J\in S(L)$ has stabilizer conjugate to $G_J$. A cube of type $J$ has dimension \(|J|\). In particular, the stabilizer of a $0$-cube is conjugate to \(G_\varnothing=\prod_{i\in I}H_i\). Applying Corollary~\ref{cor:CLL-gd-bound} to the action of $G$ on $D_G$ yields
\[
\operatorname{cat}_{\mathcal F}(G)
\leq
\max
\left\{
\operatorname{cat}_{\mathcal F|_{G_\varnothing}}
(G_\varnothing),
\;
\max_{\varnothing\neq J\in S(L)}
\left(
\operatorname{gd}(G_J)+|J|
\right)
\right\}.
\]
\end{proof}

\begin{corollary}
\label{cor:ordinary-from-pairs}
Suppose that \(H_i=\{1\}\) for every \(i\in I\). Then \(G=\mathop{\ast}_{L_1}G_i\) and \(G_J=\prod_{j\in J}G_j\). Hence
\[
\operatorname{cat}_{\mathcal F}(G)
\leq
\max_{\varnothing\neq J\in S(L)}
\left\{
\operatorname{gd}
\left(
\prod_{j\in J}G_j
\right)
+
|J|
\right\}.
\]
\end{corollary}

\begin{proof}
This is the special case \(H_i=\{1\}\) of Theorem~\ref{thm:pair-combination}. In this case $G_\varnothing$ is trivial, and hence \(\operatorname{cat}_{\mathcal F}(G_\varnothing)=0\).
\end{proof}

\section{Applications}
\label{sec:applications}

We now record several consequences of the combination estimates obtained in the previous section. We concentrate on three invariants for which the $\mathcal F$-category formalism has particularly useful geometric interpretations: amenable category, simplicial volume, and Farber's topological complexity.

\subsection{Amenable category}
\label{subsec:amenable-category}

Let $\mathrm{Am}$ denote the family of amenable subgroups of a group $G$. The corresponding invariant \(\operatorname{cat}_{\mathrm{Am}}(G)\) measures the complexity of $BG$ with respect to open covers whose fundamental groups have amenable image in $G$.

Amenable category is particularly useful because it is closely related to bounded cohomology and simplicial volume. We first specialize the combination estimates of the previous section to the family $\mathrm{Am}$.

\begin{corollary}
\label{cor:amenable-relative}
Let \(Z=Z_L(\underline A,\underline B)\) be an aspherical polyhedral product satisfying the hypotheses of Theorem~\ref{thm:polyhedral-product-combination}, and put \(G=\pi_1(Z)\). Then
\[
\operatorname{cat}_{\mathrm{Am}}(G)
\leq
\max
\left\{
\sup_{\sigma\in\Sigma_0}
\operatorname{cat}_{\mathrm{Am}|_{H_\sigma}}(H_\sigma),
\;
\sup_{\substack{\sigma\subseteq\mathcal B_L\\
\dim\sigma>0}}
\left(
\operatorname{gd}(H_\sigma)+\dim\sigma
\right)
\right\},
\]
where b\(H_\sigma=\operatorname{Stab}_G(\sigma)\).
\end{corollary}

\begin{proof}
This is Corollary~\ref{cor:polyhedral-product-gd-combination} applied to the family $\mathcal F=\mathrm{Am}$.
\end{proof}

In the transitive case, the stabilizers of the standard right-angled building are explicit, and the preceding estimate can be written entirely in terms of the defining pairs.

\begin{corollary}
\label{cor:amenable-pairs}
Let $L$ be a finite flag complex and let \(G=\Pi_{L_1}(G_i,H_i)\) be a graph product of pairs. For $J\in S(L)$, put \(G_J=\prod_{j\in J}G_j\times\prod_{i\notin J}H_i\). Then
\[
\operatorname{cat}_{\mathrm{Am}}(G)\leq\max\left\{\operatorname{cat}_{\mathrm{Am}}\left(\prod_{i\in I}H_i\right),\;\max_{\varnothing\neq J\in S(L)}\left(\operatorname{gd}(G_J)+|J|\right)\right\}.
\]
\end{corollary}

\begin{proof}
By Proposition~\ref{prop:pair-stabilizers}, the stabilizer of a $0$-cube is conjugate to \(G_\varnothing=\prod_{i\in I}H_i\), while a cube of type $J\in S(L)$ has dimension $|J|$ and stabilizer conjugate to $G_J$. The result follows from Theorem~\ref{thm:pair-combination} with $\mathcal F=\mathrm{Am}$.
\end{proof}

The estimate simplifies when the base stabilizer is amenable.

\begin{corollary}
\label{cor:amenable-pairs-amenable-H}
Under the hypotheses of Corollary~\ref{cor:amenable-pairs}, suppose that every $H_i$ is amenable. Then
\[
\operatorname{cat}_{\mathrm{Am}}(G)\leq\max_{\varnothing\neq J\in S(L)}\left(\operatorname{gd}\left(\prod_{j\in J}G_j\times\prod_{i\notin J}H_i\right)+|J|\right).
\]
\end{corollary}

\begin{proof}
Finite products of amenable groups are amenable. Hence \(G_\varnothing=\prod_{i\in I}H_i\) is amenable and therefore \(\operatorname{cat}_{\mathrm{Am}}(G_\varnothing)=0\). The result follows from Corollary~\ref{cor:amenable-pairs}.
\end{proof}

There is a different estimate when all local groups are amenable. Unlike the preceding bounds, this estimate uses the additive form of the combination theorem.

\begin{corollary}
\label{cor:all-local-amenable}
Let $L$ be a finite flag complex and let \(G=\Pi_{L_1}(G_i,H_i)\). If every $G_i$ is amenable, then \(\operatorname{cat}_{\mathrm{Am}}(G)\leq\dim L+1\).
\end{corollary}

\begin{proof}
Since $H_i\leq G_i$, every $H_i$ is amenable. Hence for every $J\in S(L)$ the local group \(G_J=\prod_{j\in J}G_j\times\prod_{i\notin J}H_i\) is amenable. Therefore \(\operatorname{cat}_{\mathrm{Am}|_{G_J}}(G_J)=0\).

Applying Corollary~\ref{cor:CLL-F-bound} to the standard right-angled building, and grouping the cells according to their dimension, gives
\[
\operatorname{cat}_{\mathrm{Am}}(G)\leq\sum_{r=1}^{\dim L+1}1=\dim L+1.
\]
\end{proof}

\begin{remark}
\label{rem:amenable-vertex-groups}
Even when all the vertex groups $G_i$ are amenable, the graph product of pairs $G$ need not itself be amenable. Thus Corollary~\ref{cor:all-local-amenable} is not a consequence of amenability of the total group. Rather, it bounds the amenable category of a potentially non-amenable group in terms of the dimension of the defining simplicial complex.
\end{remark}

The preceding results require no normality or splitting assumptions on the pairs $(G_i,H_i)$ and therefore apply to graph products of pairs in considerable generality. They provide, however, only upper bounds. We next show that for normal split pairs with amenable kernel these estimates can be replaced by exact computations.

\subsection{Amenable extensions and exact computations}
\label{subsec:amenable-extensions}

The estimates of the previous subsection apply to arbitrary graph products of pairs, but in general they provide only upper bounds. For normal pairs, the graph product of pairs admits a natural quotient which is an ordinary graph product. When the kernels are amenable, this quotient can be used to obtain substantially sharper estimates, and in the split case it leads to exact computations.

We begin with a general observation about amenable extensions.

\begin{lemma}
\label{lem:amenable-extension-upper}
Let
\[
1\longrightarrow N
\longrightarrow G
\stackrel{\rho}{\longrightarrow}
Q
\longrightarrow 1
\]
be a short exact sequence of groups. If $N$ is amenable, then \(\operatorname{cat}_{\mathrm{Am}}(G)\leq\operatorname{cat}_{\mathrm{Am}}(Q)\).
\end{lemma}

\begin{proof}
Let \(n=\operatorname{cat}_{\mathrm{Am}}(Q)\) and choose an amenable categorical cover \(BQ=U_0\cup\cdots\cup U_n\). Let \(f\colon BG\longrightarrow BQ\) be a map inducing $\rho$ on fundamental groups, and put \(V_i=f^{-1}(U_i)\). Then the sets $V_0,\ldots,V_n$ form an open cover of $BG$.

Let $V$ be a path component of some $V_i$, choose $x\in V$, and set \(\Gamma=\operatorname{im}\bigl(\pi_1(V,x)\longrightarrow G\bigr)\). Since $U_i$ is amenable categorical, the subgroup \(\rho(\Gamma)\leq Q\) is amenable. Moreover, \(\Gamma\cap N\) is amenable because it is a subgroup of the amenable group $N$. The exact sequence
\[
1\longrightarrow
\Gamma\cap N
\longrightarrow
\Gamma
\longrightarrow
\rho(\Gamma)
\longrightarrow 1
\]
shows that $\Gamma$ is amenable, since amenability is closed under extensions. Thus each $V_i$ is amenable categorical in $BG$, and hence \(\operatorname{cat}_{\mathrm{Am}}(G) \leq n\).
\end{proof}

\begin{remark}
\label{rem:extension-class}
The same argument applies more generally to a class of groups $\mathcal C$ which is closed under taking subgroups and extensions, using the family of subgroups belonging to $\mathcal C$. We shall only need the amenable case here.
\end{remark}

We now apply Lemma~\ref{lem:amenable-extension-upper} to graph products of pairs.

\begin{prop}
\label{prop:normal-pair-extension}
Let $L$ be a finite flag complex and suppose that \(H_i\trianglelefteq G_i\) for every vertex $i$. Put \(Q_i=G_i/H_i\). Then there is a natural short exact sequence
\[
1\longrightarrow
\prod_{i\in I}H_i
\longrightarrow
\Pi_{L_1}(G_i,H_i)
\stackrel{\rho}{\longrightarrow}
\mathop{\ast}_{L_1}Q_i
\longrightarrow 1.
\]
\end{prop}

\begin{proof}
Recall that the graph product of pairs is obtained from the local groups
\[
G_J=
\prod_{j\in J}G_j
\times
\prod_{i\notin J}H_i,
\qquad
J\in S(L),
\]
together with the natural inclusions associated with inclusions of simplices (see Definition~\ref{def:graph-prod-pairs}).

Since each $H_i$ is normal in $G_i$, the subgroup \(N=\prod_{i\in I}H_i\) is normal in every local group and hence in $\Pi_{L_1}(G_i,H_i)$. Moreover, \(G_J/N\cong\prod_{j\in J}Q_j\). Passing to the corresponding colimit therefore gives \(\Pi_{L_1}(G_i,H_i)/N\cong\mathop{\ast}_{L_1}Q_i\), because the colimit of the groups $\prod_{j\in J}Q_j$, for $J\in S(L)$, is precisely the ordinary graph product associated with $L_1$. This gives the claimed exact sequence.
\end{proof}

Consequently, normal amenable pairs already give an improvement over the general estimates of the previous subsection.

\begin{corollary}
\label{cor:normal-amenable-pairs}
Under the hypotheses of Proposition~\ref{prop:normal-pair-extension}, assume in addition that every $H_i$ is amenable. Then
\[
\operatorname{cat}_{\mathrm{Am}}\Pi_{L_1}(G_i,H_i)\leq\operatorname{cat}_{\mathrm{Am}}\left(\mathop{\ast}_{L_1}Q_i\right).
\]
\end{corollary}

\begin{proof}
The kernel \(\prod_{i\in I}H_i\) in Proposition~\ref{prop:normal-pair-extension} is amenable, because the vertex set is finite and finite products of amenable groups are amenable. The result follows from Lemma~\ref{lem:amenable-extension-upper}.
\end{proof}

In the split case the preceding inequality is an equality.

\begin{lemma}
\label{lem:split-pair-section}
Suppose that \(G_i=H_i\rtimes Q_i\) for every vertex $i$. Then the quotient homomorphism \(\rho\colon\Pi_{L_1}(G_i,H_i) \longrightarrow\mathop{\ast}_{L_1}Q_i\) of Proposition~\ref{prop:normal-pair-extension} admits a section. In particular, \(\mathop{\ast}_{L_1}Q_i\) embeds as a subgroup of $\Pi_{L_1}(G_i,H_i)$.
\end{lemma}

\begin{proof}
For every $i$, let \(s_i\colon Q_i\longrightarrow G_i\) be the section associated with the semidirect-product decomposition. Whenever $\{i,j\}$ is an edge of $L_1$, the groups $G_i$ and $G_j$ occur as commuting factors in the corresponding local group. Consequently, \([s_i(Q_i),s_j(Q_j)]=1\). By the universal property of the ordinary graph product, the homomorphisms $s_i$ therefore induce a homomorphism
\[
s\colon \mathop{\ast}_{L_1}Q_i \longrightarrow \Pi_{L_1}(G_i,H_i).
\]
The composite $\rho\circ s$ restricts to the identity on every vertex group $Q_i$, and hence \(\rho\circ s=\operatorname{id}\). Thus $s$ is a section and, in particular, is injective.
\end{proof}

\begin{theorem}
\label{thm:amenable-cat-split-pairs}
Let $L$ be a finite flag complex and suppose that \(G_i=H_i\rtimes Q_i\) for every vertex $i$, where every $H_i$ is amenable. Then \(\operatorname{cat}_{\mathrm{Am}} \Pi_{L_1}(G_i,H_i)=\operatorname{cat}_{\mathrm{Am}}\left(\mathop{\ast}_{L_1}Q_i\right)\).
\end{theorem}

\begin{proof}
By Proposition~\ref{prop:normal-pair-extension}, there is a short exact sequence
\[
1\longrightarrow
\prod_{i\in I}H_i
\longrightarrow
\Pi_{L_1}(G_i,H_i)
\longrightarrow
\mathop{\ast}_{L_1}Q_i
\longrightarrow1.
\]
Its kernel is amenable, and hence Lemma~\ref{lem:amenable-extension-upper} gives \(\operatorname{cat}_{\mathrm{Am}}\Pi_{L_1}(G_i,H_i)\leq\operatorname{cat}_{\mathrm{Am}}\left(\mathop{\ast}_{L_1}Q_i
\right)\).

On the other hand, Lemma~\ref{lem:split-pair-section} identifies $\mathop{\ast}_{L_1}Q_i$ with a subgroup of $\Pi_{L_1}(G_i,H_i)$. Monotonicity of amenable category under passage to subgroups therefore gives \(\operatorname{cat}_{\mathrm{Am}}\left(\mathop{\ast}_{L_1}Q_i\right)\leq\operatorname{cat}_{\mathrm{Am}}\Pi_{L_1}(G_i,H_i)\). The two inequalities prove the claim.
\end{proof}

Theorem~\ref{thm:amenable-cat-split-pairs} reduces the amenable category of a large class of graph products of pairs to that of an ordinary graph product. Exact values for the latter were obtained by Li~\cite[Theorem~3.5]{Li2023}.

\begin{theorem}[Li]
\label{thm:Li-amenable-category}
Let $Q$ be a non-trivial residually finite amenable group of type $F$, and let $Q_L$ denote the graph product with vertex group $Q$ over the finite flag complex $L$. Then
\[
\operatorname{cat}_{\mathrm{Am}}(Q_L)=\operatorname{vcd}(W_L),
\]
where $W_L$ is the right-angled Coxeter group associated with $L$.
\end{theorem}

\begin{corollary}
\label{cor:split-pairs-vcd}
Let $Q$ be a non-trivial residually finite amenable group of type $F$. Suppose that, for every vertex $i$ of $L$, \(G_i=H_i\rtimes Q\), where $H_i$ is amenable. Then \(\operatorname{cat}_{\mathrm{Am}}\Pi_{L_1}(G_i,H_i)=\operatorname{vcd}(W_L)\).
\end{corollary}

\begin{proof}
The quotient graph product is \(\mathop{\ast}_{L_1}Q=Q_L\). Hence Theorem~\ref{thm:amenable-cat-split-pairs} gives \[\operatorname{cat}_{\mathrm{Am}} \Pi_{L_1}(G_i,H_i)=\operatorname{cat}_{\mathrm{Am}}(Q_L).\] The conclusion follows from Theorem~\ref{thm:Li-amenable-category}.
\end{proof}

\begin{corollary}
\label{cor:split-pairs-RAAG}
Let $L$ be a finite flag complex and suppose that \(G_i=H_i\rtimes\mathbb Z\) for every vertex $i$, where each $H_i$ is amenable. Then
\[
\operatorname{cat}_{\mathrm{Am}}
\Pi_{L_1}(G_i,H_i)=\operatorname{vcd}(W_L).
\]
\end{corollary}

\begin{proof}
This is Corollary~\ref{cor:split-pairs-vcd} with $Q=\mathbb Z$. Equivalently, the quotient graph product is the right-angled Artin group \(A_L=\mathbb Z_L\).
\end{proof}

The splitting hypothesis in
Theorem~\ref{thm:amenable-cat-split-pairs} is used only to obtain the reverse inequality by embedding the quotient graph product into the graph product of pairs. Without a splitting, Corollary~\ref{cor:normal-amenable-pairs} still gives the upper bound, but equality is no longer automatic.

\begin{question}
\label{ques:amenable-extension-equality}
Let
\[
1\longrightarrow N\longrightarrow G\longrightarrow Q\longrightarrow1
\]
be an extension with amenable kernel. Under what additional hypotheses does one have \(\operatorname{cat}_{\mathrm{Am}}(G)=\operatorname{cat}_{\mathrm{Am}}(Q)\)? In particular, when does equality hold in Corollary~\ref{cor:normal-amenable-pairs} for non-split graph products of pairs?
\end{question}

\subsection{Join decompositions and product behaviour}
\label{subsec:join-products}

Join decompositions play a fundamental role in the topology of polyhedral products. If a simplicial complex decomposes as a join, then the corresponding polyhedral product decomposes as a Cartesian product. In the setting of moment-angle manifolds, the relationship between join decompositions and product decompositions has been studied in particular by Amelotte and Briggs \cite{AmelotteBriggs}. Here we record the corresponding behaviour for the graph products of pairs arising from polyhedral products and examine its interaction with the amenable-category computations of the previous subsection.

Let \(L=L^{(1)}*L^{(2)}\) be a join decomposition, where the vertex sets $I_1$ and $I_2$ of $L^{(1)}$ and $L^{(2)}$, respectively, are disjoint and satisfy \(I=I_1\sqcup I_2\). For $r=1,2$, write \((\underline A^{(r)},\underline B^{(r)})=\{(A_i,B_i)\}_{i\in I_r}\).

\begin{prop}
\label{prop:join-polyhedral-product}
With the notation above, there is a natural homeomorphism
\[
Z_L(\underline A,\underline B)
\cong
Z_{L^{(1)}}(\underline A^{(1)},\underline B^{(1)})
\times
Z_{L^{(2)}}(\underline A^{(2)},\underline B^{(2)}).
\]
\end{prop}

\begin{proof}
Every simplex of the join $L^{(1)}*L^{(2)}$ is uniquely of the form \(\sigma_1\cup\sigma_2\), \(\sigma_r\in L^{(r)}\). For such a simplex, the corresponding subspace in the defining union of the polyhedral product splits as
\[
\prod_{i\in\sigma_1}A_i
\times
\prod_{i\in I_1\setminus\sigma_1}B_i
\times
\prod_{i\in\sigma_2}A_i
\times
\prod_{i\in I_2\setminus\sigma_2}B_i.
\]
Taking the union over all pairs
$(\sigma_1,\sigma_2)\in L^{(1)}\times L^{(2)}$ gives
\[
Z_{L^{(1)}*L^{(2)}}(\underline A,\underline B)
=
Z_{L^{(1)}}(\underline A^{(1)},\underline B^{(1)})
\times
Z_{L^{(2)}}(\underline A^{(2)},\underline B^{(2)}).
\]
\end{proof}

The same decomposition is visible at the level of fundamental groups. In particular, in the transitive case it yields a direct product decomposition of graph products of pairs.

\begin{prop}
\label{prop:join-product-pairs}
Let \(L=L^{(1)}*L^{(2)}\) be as above, and let \(G=\Pi_{L_1}(G_i,H_i)\). For $r=1,2$, put \(G^{(r)}=\Pi_{(L^{(r)})_1}(G_i,H_i)_{i\in I_r}\). Then there is a natural isomorphism \(G\cong G^{(1)}\times G^{(2)}\).
\end{prop}

\begin{proof}
Since every vertex of $L^{(1)}$ is joined by an edge to every vertex of $L^{(2)}$, all groups associated with vertices in $I_1$ commute with all groups associated with vertices in $I_2$ in the graph product of pairs (see Definition~\ref{def:graph-prod-pairs}). Hence the canonical homomorphisms
\[
G^{(1)}\longrightarrow G,
\qquad
G^{(2)}\longrightarrow G
\]
have commuting images and induce a homomorphism \(G^{(1)}\times G^{(2)} \longrightarrow G\). Conversely, the defining local groups for $\Pi_{L_1}(G_i,H_i)$ split according to the partition $I=I_1\sqcup I_2$, and the universal property of the graph product of pairs gives a homomorphism \(G\longrightarrow G^{(1)}\times G^{(2)}\). The two homomorphisms are inverse to each other.
\end{proof}

\begin{remark}
\label{rem:join-three-levels}
Thus a join decomposition is detected simultaneously at three levels: \(L=L^{(1)}*L^{(2)}\) induces \(Z_L(\underline A,\underline B)
\cong
Z_{L^{(1)}}(\underline A^{(1)},\underline B^{(1)})
\times
Z_{L^{(2)}}(\underline A^{(2)},\underline B^{(2)})\) and, in the transitive case, \(\Pi_{L_1}(G_i,H_i)\cong G^{(1)}\times G^{(2)}\). This is the group-theoretic counterpart of the familiar join--product correspondence for polyhedral products.
\end{remark}

We next combine Proposition~\ref{prop:join-product-pairs} with the exact computation from Theorem~\ref{thm:amenable-cat-split-pairs}. The resulting behaviour is more subtle than ordinary additivity.

\begin{corollary}
\label{cor:join-split-amenable}
Suppose that \(L=L^{(1)}*L^{(2)}\) and that \(G_i=H_i\rtimes Q_i\) for every vertex $i$, with each $H_i$ amenable. Then
\[
\operatorname{cat}_{\mathrm{Am}}
\Pi_{L_1}(G_i,H_i)
=
\operatorname{cat}_{\mathrm{Am}}
\left(
\mathop{\ast}_{L_1} Q_i
\right).
\]
Moreover,
\[
\mathop{\ast}_{L_1}Q_i
\cong
\left(
\mathop{\ast}_{(L^{(1)})_1}Q_i
\right)
\times
\left(
\mathop{\ast}_{(L^{(2)})_1}Q_i
\right).
\]
\end{corollary}

\begin{proof}
The first assertion is Theorem~\ref{thm:amenable-cat-split-pairs}. For the second, every vertex of $L^{(1)}$ is adjacent to every vertex of $L^{(2)}$. Hence the defining relations of the ordinary graph product force the two corresponding graph-product factors to commute, giving the asserted direct-product decomposition.
\end{proof}

When all quotient groups coincide with a group to which Li's theorem applies, the preceding statement gives an exact combinatorial description.

\begin{corollary}
\label{cor:join-vcd}
Let $Q$ be a non-trivial residually finite amenable group of type $F$, and suppose that \(G_i=H_i\rtimes Q\) for every vertex $i$, where each $H_i$ is amenable. If \(L=L^{(1)}*L^{(2)}\), then \(\operatorname{cat}_{\mathrm{Am}} \Pi_{L_1}(G_i,H_i)=\operatorname{vcd}(W_L)\). Furthermore, \(W_L\cong W_{L^{(1)}}\times W_{L^{(2)}}\). Consequently,
\[
\operatorname{cat}_{\mathrm{Am}}
\Pi_{L_1}(G_i,H_i)
=
\operatorname{vcd}
\left(
W_{L^{(1)}}\times W_{L^{(2)}}
\right).
\]
\end{corollary}

\begin{proof}
By Corollary~\ref{cor:split-pairs-vcd},
\[
\operatorname{cat}_{\mathrm{Am}}\Pi_{L_1}(G_i,H_i)=\operatorname{vcd}(W_L).
\]
The defining graph of the right-angled Coxeter group associated with a join is itself a join, and therefore \(W_{L^{(1)}*L^{(2)}} \cong W_{L^{(1)}}\times W_{L^{(2)}}\). This proves the claim.
\end{proof}

\begin{remark}
\label{rem:join-nonadditivity}
It is important that the corollary \ref{cor:join-vcd} is formulated in terms of \(\operatorname{vcd}\left(W_{L^{(1)}}\times W_{L^{(2)}}\right)\) rather than as an additivity statement. Amenable category is not additive under direct products in general. Indeed, Li \cite{Li2023} exhibits right-angled Artin groups for which
\[
\operatorname{cat}_{\mathrm{Am}}(A_{L^{(1)}}\times A_{L^{(2)}})<\operatorname{cat}_{\mathrm{Am}}(A_{L^{(1)}})+\operatorname{cat}_{\mathrm{Am}}(A_{L^{(2)}}).
\]
Thus join decompositions of the defining simplicial complex give canonical product decompositions of both the polyhedral product and its graph-product-of-pairs fundamental group, while the numerical behaviour of amenable category under these products can be genuinely non-additive.
\end{remark}

The preceding discussion mirrors a familiar phenomenon in toric topology: combinatorial join decompositions of the defining simplicial complex are reflected geometrically by Cartesian product decompositions of the associated polyhedral product. The results above show that, for the aspherical polyhedral products considered here, the same decomposition is also visible at the level of the relative fundamental group. In contrast with the geometric product decomposition, however, amenable category need not decompose additively. This provides a group-theoretic complement to the join--product viewpoint appearing in the study of product decompositions of moment-angle manifolds.

\subsection{Simplicial volume and manifold applications}
\label{subsec:simplicial-volume}

We now turn the amenable-category estimates obtained above into
vanishing results for simplicial volume. An important point is that, at this stage, the manifold itself need not be aspherical: it is enough to control the amenable category of its fundamental group.

Recall that if $M$ is an oriented closed connected $m$-manifold, its \emph{simplicial volume} $\|M\|$ is the infimum of the $\ell^1$-norms of real singular cycles representing the fundamental class $[M]\in H_m(M;\mathbb R)$.

We shall use the following classical vanishing theorem.

\begin{theorem}[Gromov {\cite{gromov1982volume}}]
\label{thm:Gromov-amenable-vanishing}
Let $M$ be an oriented closed connected manifold. If \(\operatorname{cat}_{\mathrm{Am}}(M)<\dim M\), then \(\|M\|=0\).
\end{theorem}

\begin{lemma}
\label{lem:space-group-amenable-category}
Let $X$ be a path-connected CW-complex. Then \(\operatorname{cat}_{\mathrm{Am}}(X)\leq\operatorname{cat}_{\mathrm{Am}}(\pi_1(X))\).
\end{lemma}

\begin{proof}
Let \(c_X\colon X\longrightarrow B\pi_1(X)\) be a classifying map inducing the identity on fundamental groups, and let \(B\pi_1(X)=U_0\cup\cdots\cup U_n\) be an amenable categorical cover, where \(n=\operatorname{cat}_{\mathrm{Am}}(\pi_1(X))\). Put \(V_i=c_X^{-1}(U_i)\). Then the $V_i$ form an open cover of $X$.

For every path component $V$ of some $V_i$, the image of \(\pi_1(V)\longrightarrow\pi_1(X)\) is contained in the image of the fundamental group of the corresponding component of $U_i$ in $\pi_1(X)$. The latter is amenable. Hence $V_i$ is amenable categorical in $X$ for every $i$, and therefore \(\operatorname{cat}_{\mathrm{Am}}(X)\leq n\).
\end{proof}

\begin{prop}
\label{prop:group-category-simplicial-volume}
Let $M$ be an oriented closed connected $m$-manifold. If \(\operatorname{cat}_{\mathrm{Am}}(\pi_1(M))<m\), then \(\|M\|=0\).
\end{prop}

\begin{proof}
By Lemma~\ref{lem:space-group-amenable-category},
\[
\operatorname{cat}_{\mathrm{Am}}(M)\leq\operatorname{cat}_{\mathrm{Am}}(\pi_1(M))<m.
\]
The conclusion follows from Theorem~\ref{thm:Gromov-amenable-vanishing}.
\end{proof}

\begin{theorem}
\label{thm:pairs-simplicial-volume}
Let $M$ be an oriented closed connected $m$-manifold such that \(\pi_1(M)\cong \Pi_{L_1}(G_i,H_i)\), where $L$ is a finite flag complex. For $J\in S(L)$, put \(G_J=\prod_{j\in J}G_j\times\prod_{i\notin J}H_i\). If
\[
\max
\left\{
\operatorname{cat}_{\mathrm{Am}}
\left(
\prod_{i\in I}H_i
\right),
\;
\max_{\varnothing\neq J\in S(L)}
\left(
\operatorname{gd}(G_J)+|J|
\right)
\right\}
<m,
\]
then \(\|M\|=0\).
\end{theorem}

\begin{proof}
Corollary~\ref{cor:amenable-pairs} gives
\[
\operatorname{cat}_{\mathrm{Am}}(\pi_1(M))
\leq
\max
\left\{
\operatorname{cat}_{\mathrm{Am}}
\left(
\prod_{i\in I}H_i
\right),
\;
\max_{\varnothing\neq J\in S(L)}
\left(
\operatorname{gd}(G_J)+|J|
\right)
\right\}.
\]
The assumed strict inequality and Proposition~\ref{prop:group-category-simplicial-volume} give the claim.
\end{proof}

\begin{corollary}
\label{cor:amenable-local-simplicial-volume}
Let $M$ be an oriented closed connected $m$-manifold such that \(\pi_1(M)\cong\Pi_{L_1}(G_i,H_i)\), where $L$ is finite and flag. If every $G_i$ is amenable and \(\dim L+1<m\), then \(\|M\|=0\).
\end{corollary}

\begin{proof}
By Corollary~\ref{cor:all-local-amenable}, \(\operatorname{cat}_{\mathrm{Am}}(\pi_1(M))\leq\dim L+1<m\). Apply Proposition~\ref{prop:group-category-simplicial-volume}.
\end{proof}

The exact computations of Subsection~\ref{subsec:amenable-extensions} yield a sharper criterion for manifolds whose fundamental groups are split graph products of pairs.

\begin{theorem}
\label{thm:split-pairs-simplicial-volume}
Let $M$ be an oriented closed connected $m$-manifold and suppose that \(\pi_1(M)\cong\Pi_{L_1}(H_i\rtimes Q_i,H_i)\), where every $H_i$ is amenable. If \(\operatorname{cat}_{\mathrm{Am}}
\left(\mathop{\ast}_{L_1}Q_i \right)<m\), then \(\|M\|=0\).
\end{theorem}

\begin{proof}
By Theorem~\ref{thm:amenable-cat-split-pairs},
\[
\operatorname{cat}_{\mathrm{Am}}(\pi_1(M))
=
\operatorname{cat}_{\mathrm{Am}}
\left(
\mathop{\ast}_{L_1}Q_i
\right).
\]
The result follows from Proposition~\ref{prop:group-category-simplicial-volume}.
\end{proof}

\begin{corollary}
\label{cor:split-pairs-vcd-simplicial-volume}
Let $Q$ be a non-trivial residually finite amenable group of type $F$, and let $M$ be an oriented closed connected $m$-manifold with \(\pi_1(M)\cong\Pi_{L_1}(H_i\rtimes Q,H_i)\), where every $H_i$ is amenable. If \(\operatorname{vcd}(W_L)<m\), then \(\|M\|=0\).
\end{corollary}

\begin{proof}
By Corollary~\ref{cor:split-pairs-vcd}, \(\operatorname{cat}_{\mathrm{Am}}(\pi_1(M))=\operatorname{vcd}(W_L)\). Hence the assumption implies \[\operatorname{cat}_{\mathrm{Am}}(\pi_1(M))<m,\] and Proposition~\ref{prop:group-category-simplicial-volume} applies.
\end{proof}

\begin{corollary}
\label{cor:RAAG-pairs-simplicial-volume}
Let $M$ be an oriented closed connected $m$-manifold with \(\pi_1(M)\cong\Pi_{L_1}(H_i\rtimes\mathbb Z,H_i)\), where every $H_i$ is amenable. If \(\operatorname{vcd}(W_L)<m\), then \(\|M\|=0\).
\end{corollary}

\begin{remark}
\label{rem:simplicial-volume-symmetry}
The group-theoretic vanishing criteria above should be distinguished from the classical vanishing results coming from continuous group actions. For example, whenever a complex moment-angle complex \(\mathcal Z_L=Z_L(D^2,S^1)\) is a closed manifold, it carries the canonical coordinatewise torus action. Thus its simplicial-volume behaviour is already constrained by the classical vanishing results for manifolds with non-trivial circle actions.

The real moment-angle case is different. A real moment-angle manifold \(\mathcal Z_L^{\mathbb R}=Z_L(D^1,S^0)\) carries the canonical action of the finite group $(\mathbb Z/2)^m$, rather than a torus action, so circle-action arguments do not apply. Nevertheless, if $L$ is a triangulated sphere, then \(\dim \mathcal Z_L^{\mathbb R}=\dim L+1\), and hence Corollary~\ref{cor:amenable-local-simplicial-volume} does not yield vanishing, since its hypothesis requires the strict inequality \(\dim L+1<\dim M\). Thus the criterion is genuinely selective rather than a general vanishing statement for moment-angle manifolds.
\end{remark}

\begin{remark}
\label{rem:Davis-simplicial-volume}
A natural testing ground for the preceding criteria is provided by manifolds obtained from Davis' reflection-group construction. Their simplicial volume is not uniformly zero, and hence no blanket vanishing statement can hold for this class. It would therefore be interesting to determine whether the amenable-category estimates above single out subclasses of Davis-type manifolds for which the simplicial volume vanishes.
\end{remark}

\section{Polyhedral products of classifying spaces}
\label{sec:classifying-spaces}

The preceding results become particularly concrete for polyhedral products built from classifying spaces of pairs of groups. Let $G_i$ be a discrete group and let \(H_i\leq G_i\), \(i\in I\). Choose CW-models for $BH_i$ and $BG_i$ such that the homomorphism induced by the inclusion $H_i\hookrightarrow G_i$ is represented by a cellular inclusion \(BH_i\hookrightarrow BG_i\). If necessary, this can always be achieved by replacing the map $BH_i\to BG_i$ by its mapping cylinder. We consider the polyhedral product \(Z_L(\underline{BG},\underline{BH})=Z_L\bigl((BG_i,BH_i)_{i\in I}\bigr)\).

The advantage of this class is that both the fundamental-group data and the stabilizers of the associated building admit purely group-theoretic descriptions.

\subsection{The pairs $(BG_i,BH_i)$}
\label{subsec:BG-BH}

Let \(p_i:EG_i\longrightarrow BG_i\) be the universal covering. The inverse image \(p_i^{-1}(BH_i)\) is a disjoint union of contractible components. The group $G_i$ acts transitively on the set of these components, and the stabilizer of a chosen component is $H_i$. Consequently, \(\pi_0\bigl(p_i^{-1}(BH_i)\bigr) \cong G_i/H_i\) as a $G_i$-set.

Thus the fundamental-group data associated with the pair $(BG_i,BH_i)$ are precisely \((G_i,G_i/H_i)\). In particular, the local permutation action is transitive, and hence this class of polyhedral products lies exactly in the graph-product-of-pairs situation considered in Subsection~\ref{subsec:combination-pairs}.

We collect the basic consequences.

\begin{prop}
\label{prop:BG-BH-basic}
Let $L$ be a simplicial complex with finite vertex set $I$, and let $H_i\leq G_i$ for every $i\in I$. Then:
\begin{enumerate}
    \item the fundamental-group data associated with
    $(BG_i,BH_i)$ are \(E_i\cong G_i/H_i\);

    \item there is a natural isomorphism \(\pi_1 Z_L(\underline{BG},\underline{BH})\cong\Pi_{L_1}(G_i,H_i)\);

    \item if $H_i<G_i$ is proper for every $i\in I$, then \(Z_L(\underline{BG},\underline{BH})\) is aspherical if and only if $L$ is flag.
\end{enumerate}
\end{prop}

\begin{proof}
The universal cover of $BG_i$ is $EG_i$. Since $BH_i\to BG_i$ induces the inclusion \(H_i\hookrightarrow G_i\) on fundamental groups, every component of $p_i^{-1}(BH_i)$ is a universal cover of $BH_i$ and is therefore contractible.

The deck-transformation action of $G_i$ permutes these components transitively. Choosing one component, its stabilizer is $H_i$, and hence \(E_i=\pi_0\bigl(p_i^{-1}(BH_i)\bigr)\cong G_i/H_i\). This proves~(1).

By Theorem~\ref{thm:Davis-fundamental-group}, \(\pi_1 Z_L(\underline{BG},\underline{BH}) \cong Z_{L_1}(\underline G,\underline E)\). Since each action \(G_i\curvearrowright G_i/H_i\) is transitive, Davis's identification of transitive relative graph products with the graph products of pairs of Januszkiewicz--Świątkowski gives \(Z_{L_1}(\underline G,\underline E)\cong\Pi_{L_1}(G_i,H_i)\), which proves~(2).

Finally, suppose that every $H_i$ is proper in $G_i$. Then \(|G_i/H_i|>1\) for every $i$, and hence all vertices contribute to the nontrivial right-angled structure. Moreover, both $BG_i$ and $BH_i$ are aspherical and each inclusion \(BH_i\longrightarrow BG_i\) is $\pi_1$-injective. Thus every pair $(BG_i,BH_i)$ is an aspherical pair. Corollary~\ref{cor:asphericity-standard} therefore shows that the only remaining condition for asphericity of the polyhedral product is the flagness of $L$.
\end{proof}

\begin{remark}
The properness assumption in Proposition~\ref{prop:BG-BH-basic}(3) is used only to guarantee that \(|G_i/H_i|>1\) for every $i$. If $H_i=G_i$ for some indices, then the corresponding $G_i$-sets consist of a single point. In this case one must pass, as in Subsection~\ref{subsec:asphericity}, to the reduced vertex set \(I''=\{i\in I\mid H_i\neq G_i\}\) and the corresponding full subcomplex \(L''=L_{I''}\). The relevant right-angled building is then indexed by $L''$.
\end{remark}

When $L$ is flag and all $H_i$ are proper, the universal covering of the polyhedral product can therefore be studied through the right-angled building associated with the graph product of pairs \(\Pi_{L_1}(G_i,H_i)\). For every $J\in S(L)$, put \(G_J=\prod_{j\in J}G_j\times\prod_{i\in I\setminus J}H_i\). By Proposition~\ref{prop:pair-stabilizers}, the stabilizer of a cube of type $J$ is conjugate to $G_J$. In particular, \(G_\varnothing=\prod_{i\in I}H_i\).

We therefore obtain directly from Theorem~\ref{thm:pair-combination} the following form of the main estimate.

\begin{corollary}
\label{cor:BG-BH-combination}
Let $L$ be a finite flag complex and let \(H_i<G_i\) for every $i\in I$. Put \(Z=Z_L(\underline{BG},\underline{BH})\). Then, for every family $\mathcal F$ of subgroups of \(G=\pi_1(Z) \), one has
\[
\operatorname{cat}_{\mathcal F}(G)\leq\max\left\{\operatorname{cat}_{\mathcal F|_{G_\varnothing}}
(G_\varnothing),\;\max_{\varnothing\neq J\in S(L)}
\left(\operatorname{gd}(G_J)+|J|\right)\right\},
\]
where \(G_\varnothing=\prod_{i\in I}H_i\) and
\(G_J=\prod_{j\in J}G_j
\times\prod_{i\notin J}H_i\).
\end{corollary}

\begin{proof}
By Proposition~\ref{prop:BG-BH-basic}, \(G\cong\Pi_{L_1}(G_i,H_i)\), and the result is therefore precisely Theorem~\ref{thm:pair-combination}.
\end{proof}

\begin{example}
\label{ex:BG-point}
If $H_i=\{1\}$ for every $i$, then \(BH_i=*\) and \(Z_L(\underline{BG},*)\) is aspherical precisely when $L$ is flag. Its fundamental group is the ordinary graph product \(\pi_1Z_L(\underline{BG},*)\cong\mathop{\ast}_{L_1}G_i\). The local groups reduce to \(G_J=\prod_{j\in J}G_j\), recovering the ordinary graph-product case.
\end{example}

\begin{example}
\label{ex:RAAG-classifying}
Taking \(G_i=\mathbb Z\), \(H_i=\{1\}\), gives \(Z_L(S^1,*)\simeq BA_{L_1}\) when $L$ is flag, where $A_{L_1}$ denotes the associated right-angled Artin group.
\end{example}

\begin{example}[Real Davis--Januszkiewicz space]
\label{ex:RACG-classifying}
Taking \(G_i=\mathbb Z/2,\qquad H_i=\{1\}\), for every $i\in I$, gives \(Z_L(\mathbb RP^\infty,*)\). If $L$ is flag, Proposition~\ref{prop:BG-BH-basic} implies that this space is aspherical, and
\[
\pi_1 Z_L(\mathbb RP^\infty,*)\cong\mathop{\ast}_{L_1}\mathbb Z/2=W_{L_1},
\]
where $W_{L_1}$ is the right-angled Coxeter group associated with the $1$-skeleton of $L$. Thus \(Z_L(\mathbb RP^\infty,*) \simeq BW_{L_1}\). This is the real Davis--Januszkiewicz space associated with $L$.
\end{example}

\section{Farrell--Jones inheritance and rigidity consequences}
\label{sec:rigidity}

The preceding sections describe the fundamental groups of polyhedral products in terms of relative graph products and graph products of pairs, and provide contractible right-angled buildings on which these groups act. We now investigate the consequences of this structure for the Farrell--Jones conjecture and for topological rigidity.

The relevant framework is the Farrell--Jones conjecture in algebraic $K$- and $L$-theory. For closed aspherical manifolds, the conjecture provides the algebraic input required in surgery theory to deduce the Borel conjecture.

Our purpose in this section is threefold. First, we recall the Farrell--Jones theorem of Gandini--R\"uping for ordinary graph products \cite{GandiniRueping}. Second, we develop a pair-deletion argument for graph products of pairs, reducing the global Farrell--Jones problem to a local combination problem associated with the inclusions \(H_i\leq G_i\). Third, we establish concrete classes of such local pairs using acylindrical Bass--Serre actions, elementary amenable extension arguments, and known global descriptions of finite group pairs. We then apply these results to fundamental groups of polyhedral products and to Borel rigidity for closed aspherical polyhedral-product manifolds.

\subsection{Manifold polyhedral products}
\label{subsec:manifold-polyhedral-products}

We shall call a polyhedral product \(M=Z_L(\underline A,\underline B)\) a \emph{manifold polyhedral product} if it is a topological manifold. In this section we are primarily interested in the case in which $M$ is closed, connected and aspherical.

If \(M=Z_L(\underline A,\underline B)\) is aspherical, then \(M\simeq B\pi_1(M)\). Consequently, the homotopy type of $M$ is determined by its fundamental group. By Theorem~\ref{thm:Davis-fundamental-group}, this group has the description
\[
\pi_1(M)
\cong
Z_{L_1}(\underline G,\underline E),
\]
where
\[
G_i=\pi_1(A(i)), \qquad E_i=\pi_0(\widetilde B(i)).
\]

Thus the rigidity problem for an aspherical manifold polyhedral product naturally becomes a rigidity problem for a relative graph product.

\begin{remark}
\label{rem:torsion-free-manifold-group}
If $M$ is a closed aspherical manifold, then \(G=\pi_1(M)\) is torsion-free and admits a finite classifying space, namely $M$ itself. In particular, \(\operatorname{cd}(G)=\dim M\).
\end{remark}

The building description developed in Subsection~\ref{subsec:building-action} gives an additional geometric structure. Under the flag hypotheses considered above, \(G\curvearrowright \mathcal B_L\) acts cellularly on a finite-dimensional contractible $\operatorname{CAT}(0)$ cubical complex. It is important, however, to distinguish the existence of a $\operatorname{CAT}(0)$ space from the existence of a geometric $\operatorname{CAT}(0)$ action.

\begin{remark}
\label{rem:not-automatically-CAT0-group}
The action \(G\curvearrowright\mathcal B_L\) need not be proper. Indeed, in the graph-product-of-pairs case the stabilizer of a cube of type $J$ is conjugate to
\[
G_J=\prod_{j\in J}G_j\times \prod_{i\notin J}H_i,
\]
which can be infinite. Hence the fact that $\mathcal B_L$ is $\operatorname{CAT}(0)$ does not by itself imply that $G$ is a $\operatorname{CAT}(0)$ group.
\end{remark}

This distinction becomes essential when applying the Farrell--Jones conjecture.

A basic example is provided by real moment-angle manifolds. If $L$ is a flag triangulation of a sphere, then \(\mathcal Z_L^{\mathbb R}=Z_L(D^1,S^0)\) is a closed aspherical manifold. Its universal cover is the Davis complex associated with the right-angled Coxeter group $W_L$, and its fundamental group is a torsion-free finite-index subgroup of $W_L$. This example acts as a prototype and illustrates the general philosophy: the right-angled geometry of the universal cover places the fundamental group within a class for which rigidity techniques become available. For further details see \cite{gkeneralis-real}. 

\subsection{Farrell--Jones}
\label{subsec:Farrell-Jones}

We briefly recall the form of the Farrell--Jones conjecture needed below. Let $\mathcal{VC}yc$ denote the family of virtually cyclic subgroups of a group $G$. The conjecture predicts that the assembly maps
\[
H_n^G \bigl( E_{\mathcal{VC}yc}G;\mathbf K_{\mathcal A}\bigr)\longrightarrow K_n(\mathcal A\rtimes G)
\]
and
\[
H_n^G\bigl( E_{\mathcal{VC}yc}G;\mathbf L_{\mathcal A}^{\langle-\infty\rangle}\bigr)\longrightarrow L_n^{\langle-\infty\rangle}
(\mathcal A\rtimes G)
\]
are isomorphisms, for suitable additive categories $\mathcal A$ with $G$-action. More precisely, throughout this subsection, unless explicitly stated otherwise, the Farrell--Jones conjecture means the \(K\)- and \(L\)-theoretic conjecture with coefficients in additive categories and with finite wreath products.

The relationship between graph products and assembly maps in algebraic $K$- and $L$-theory has been studied in considerable detail. Gandini and R\"uping proved that the class of groups satisfying the Full Farrell--Jones Conjecture is closed under graph products \cite{GandiniRueping}. Subsequently, Kasprowski, Li and L\"uck developed explicit formulas for equivariant homology theories, group homology, algebraic $K$- and $L$-theory, and topological $K$-theory of graph products \cite{KasprowskiLiLueck}. Their approach expresses the global invariants of a graph product in terms of the subgroups associated with simplices of the flag complex.

This local-to-global structure is closely related to the building description used in the present paper: for a simplex $J\in S(L)$, the stabilizer of the corresponding cube in the ordinary graph-product case is \(G_J=\prod_{j\in J}G_j\). Thus the simplex subgroups which enter the $K$- and $L$-theoretic decomposition of Kasprowski--Li--L\"uck are precisely the isotropy groups appearing in the standard right-angled building.

We shall use the version with coefficients in additive categories and finite wreath products. Following Gandini--R\"uping \cite{GandiniRueping}, let $\mathcal C_{\mathrm{FJ}}$ denote the class of groups satisfying this version of the $K$- and $L$-theoretic Farrell--Jones conjecture.

The class $\mathcal C_{\mathrm{FJ}}$ enjoys strong inheritance properties. In particular, it contains hyperbolic groups, $\operatorname{CAT}(0)$ groups and virtually polycyclic groups, and is closed under subgroups, finite-index overgroups, directed colimits, finite direct products and free products. We shall also use the following inheritance principle. Let \(f\colon G\longrightarrow Q\) be a homomorphism. If \(Q\in\mathcal C_{\mathrm{FJ}}\), \(\ker(f)\in\mathcal C_{\mathrm{FJ}}\), and \(f^{-1}(C)\in\mathcal C_{\mathrm{FJ}}\) for every infinite cyclic subgroup \(C\leq Q\), then \(G\in\mathcal C_{\mathrm{FJ}}\); see~\cite[Proposition~1.1]{GandiniRueping}.

The main result of Gandini--R\"uping relevant to the present paper is the following.

\begin{theorem}[Gandini--R\"uping {\cite[Theorem~1]{GandiniRueping}}]
\label{thm:GR-graph-products}
Let $\Gamma$ be a simplicial graph and let $(G_i)_{i\in I}$ be a family of groups satisfying the Farrell--Jones conjecture. Then the graph product \(\mathop{\ast}_{\Gamma}G_i\) also satisfies the Farrell--Jones conjecture.
\end{theorem}

Theorem~\ref{thm:GR-graph-products} gives an immediate Farrell--Jones result for those aspherical polyhedral products whose fundamental groups are ordinary graph products. We specialize this theorem to our setting.

\begin{prop}
\label{prop:FJ-ordinary-polyhedral-product}
Let \(Z=Z_L(\underline A,\underline B)\) be a polyhedral product such that the associated actions \(G_i\curvearrowright E_i\) are simply transitive. Then \(\pi_1(Z)\cong\mathop{\ast}_{L_1}G_i\). If every $G_i$ satisfies the Farrell--Jones conjecture, then $\pi_1(Z)$ satisfies the Farrell--Jones conjecture.
\end{prop}

\begin{proof}
By Proposition~\ref{prop:ordinary-graph-product-stabilizers}, the simple transitivity hypothesis identifies the relative graph product with the ordinary graph product \(\pi_1(Z)\cong\mathop{\ast}_{L_1}G_i\). The conclusion follows from Theorem~\ref{thm:GR-graph-products}.
\end{proof}

For the classifying-space case, this gives:

\begin{corollary}
\label{cor:FJ-BG-point}
Let \(L\) be a finite simplicial complex and let \((G_i)_{i\in I}\) be groups in \(\mathcal C_{\mathrm{FJ}}\). Then \(\pi_1 Z_L(\underline{BG},*)\cong\mathop{\ast}_{L_1}G_i\) belongs to \(\mathcal C_{\mathrm{FJ}}\).
\end{corollary}

The situation is subtler for graph products of pairs. If \(E_i\cong G_i/H_i\), then \(\pi_1(Z)\cong\Pi_{L_1}(G_i,H_i)\), and Theorem~\ref{thm:GR-graph-products} does not apply directly.

\subsubsection{Pair deletion}

Let $\Gamma$ be a finite simplicial graph and write \(P_\Gamma=\Pi_\Gamma(G_i,H_i)\). Fix a vertex $v\in V(\Gamma)$ and put
\[
\Gamma'=\Gamma-\{v\},
\qquad
P'= \Pi_{\Gamma'}(G_i,H_i).
\]

Recall that \(\operatorname{lk}_\Gamma(v)=\{i\neq v:\{i,v\}\in E(\Gamma)\}\) and \(\operatorname{st}_\Gamma(v)=\{v\}\cup\operatorname{lk}_\Gamma(v)\).

\begin{definition}
\label{def:relative-link-group}
The \emph{relative link subgroup} of $v$ is
\[
R_v=\left\langle G_i\;(i\in\operatorname{lk}_\Gamma(v)),\;H_i\;(i\notin\operatorname{st}_\Gamma(v))\right\rangle\leq P'.
\]
\end{definition}

When all $H_i$ are trivial, $R_v$ is the ordinary graph product associated with the link of $v$.

\begin{prop}[Pair-deletion decomposition]
\label{prop:pair-deletion-amalgam}
There is a natural isomorphism
\[
P_\Gamma
\cong
(H_v\times P')
*_{H_v\times R_v}
(G_v\times R_v).
\]
\end{prop}

\begin{proof}
Recall that the graph product of pairs is the colimit of the system of local groups \(G_J=\prod_{j\in J}G_j\times\prod_{i\notin J}H_i\) , where $J$ ranges over the cliques of $\Gamma$ (see Definition~\ref{def:graph-prod-pairs}). Every simplex containing $v$ is of the form \(\{v\}\cup\sigma\), \(\sigma\subseteq\operatorname{lk}_\Gamma(v)\). Hence the part of the colimit involving $v$ is generated by $G_v$ together with $R_v$.

By construction, $G_v$ commutes with $R_v$, while $H_v$ commutes with all of $P'$. Therefore the two natural subgroups
\[
A=H_v\times P',
\qquad
B=G_v\times R_v
\]
have common subgroup \(E=H_v\times R_v\). The defining relations of $P_\Gamma$ are precisely those of $A$ and $B$ together with the identification of $E$. The universal property of amalgamated free products gives \(P_\Gamma\cong A*_E B\).
\end{proof}

The preceding decomposition defines a canonical epimorphism \(\rho_v:P_\Gamma\longrightarrow P'\) by
\[
\rho_v(h,p)=p
\qquad
(h,p)\in H_v\times P',
\]
and
\[
\rho_v(g,r)=r
\qquad
(g,r)\in G_v\times R_v.
\]
The canonical inclusion $P'\hookrightarrow P_\Gamma$ is a section of $\rho_v$.

\begin{prop}[Kernel of pair deletion]
\label{prop:pair-deletion-kernel}
Let \(K_v=\ker(\rho_v)\). Then \(K_v\cong\underset{P'/R_v}{\ast_{H_v}}G_v\). In other words, $K_v$ is the amalgamated free product of copies of $G_v$, indexed by the left cosets $P'/R_v$, with all copies of $H_v$ identified.
\end{prop}

\begin{proof}
Let $T$ be the Bass--Serre tree associated with
\[
P_\Gamma
=
(H_v\times P')
*_{H_v\times R_v}
(G_v\times R_v).
\]
Since $K_v$ is normal and \(P_\Gamma/K_v\cong P'\), the quotient graph $K_v\backslash T$ has one vertex of $H_v\times P'$--type, while the vertices of $G_v\times R_v$--type are indexed by \(P'/R_v\). The edges are indexed by the same coset space.

The corresponding stabilizers are \(K_v\cap(H_v\times P')=H_v\), \(K_v\cap(G_v\times R_v)=G_v\), and \(K_v\cap(H_v\times R_v)=H_v\). Since $H_v$ commutes with $P'$, all edge inclusions identify the same copy of $H_v$. Bass--Serre theory therefore yields \(K_v\cong \underset{P'/R_v}{\ast_{H_v}}G_v\).
\end{proof}

\begin{prop}[Cyclic preimages]
\label{prop:pair-deletion-cyclic}
Let $C\leq P'$ be an infinite cyclic subgroup. Then \(\rho_v^{-1}(C)\cong\left(\underset{P'/R_v}{\ast_{H_v}}G_v\right)\rtimes C\), where $C$ acts on the index set $P'/R_v$ by left multiplication. Moreover, $C$ fixes the common subgroup $H_v$ pointwise.
\end{prop}

\begin{proof}
Since $\rho_v$ admits the canonical section \(P'\hookrightarrow P_\Gamma\), we have \(P_\Gamma\cong K_v\rtimes P'\). Restricting to $C\leq P'$ gives \(\rho_v^{-1}(C)\cong K_v\rtimes C\). By Proposition~\ref{prop:pair-deletion-kernel}, \(K_v\cong\underset{P'/R_v}{\ast_{H_v}}G_v\). Conjugation by $C$ permutes the factors through the left action of $C$ on $P'/R_v$.  Finally, $H_v$ commutes with $P'$, and hence the action of $C$ fixes $H_v$ pointwise.
\end{proof}

\subsubsection{FJ-admissible pairs}

The preceding propositions reduce the Farrell--Jones problem for graph products of pairs to a local combination problem. For a nonempty set \(S\), we write \(\underset{s\in S}{\ast_H}G\) for the free product of copies \(G_s\cong G\), indexed by \(S\), with all copies of the subgroup \(H\) identified by the identity map. If \(S\) is a \(\mathbb Z\)-set, the action of \(\mathbb Z\) on \(S\) induces an action on this amalgamated free product by permuting the factors and fixing the common subgroup \(H\) pointwise. We write \(\mathcal G(G,H;S)=\left(\underset{s\in S}{\ast_H}G\right)\rtimes\mathbb Z\).

\begin{definition}
\label{def:FJ-admissible-pair}
A pair \(H\leq G\) is called \emph{FJ-admissible} if \(G\in\mathcal C_{\mathrm{FJ}}\) and \(\mathcal G(G,H;S)\in\mathcal C_{\mathrm{FJ}}\) for every nonempty \(\mathbb Z\)-set \(S\).
\end{definition}

For a nonempty $\mathbb Z$--set $S$, put \(\mathcal G(G,H;S)=\left( \underset{s\in S}{\ast_H}G \right)\rtimes\mathbb Z\). The proof of \cite[Lemma~2.8]{GandiniRueping} suggests decomposing $S$ into $\mathbb Z$--orbits.  The same decomposition is available in the relative setting, but disjoint orbits are now glued along $H\times\mathbb Z$ rather than along the infinite cyclic factor alone.

\begin{lemma}
\label{lem:relative-orbit-decomposition}
Let $S=S_1\sqcup S_2$ be a decomposition into $\mathbb Z$--invariant subsets. Then
\[
\mathcal G(G,H;S)
\cong
\mathcal G(G,H;S_1)
*_{H\times\mathbb Z}
\mathcal G(G,H;S_2).
\]
Moreover:
\begin{enumerate}
\item if $S\cong\mathbb Z/m\mathbb Z$, then \(\mathcal G(G,H;S)\) contains \(\left(\underset{m}{\ast_H}G\right)\times\mathbb Z\) as a subgroup of index $m$;

\item if $S\cong\mathbb Z$, then \(\mathcal G(G,H;S)\cong G*_H(H\times\mathbb Z)\).
\end{enumerate}
\end{lemma}
\begin{proof}
Write \(K(S)=\underset{s\in S}{\ast_H}G\). If \(S=S_1\sqcup S_2\) is a decomposition into \(\mathbb Z\)-invariant subsets, then \(K(S) \cong K(S_1)*_H K(S_2)\). The generator \(t\) of \(\mathbb Z\) preserves both factors and fixes the common subgroup \(H\) pointwise. Hence
\[
K(S)\rtimes\langle t\rangle
\cong
\bigl(K(S_1)\rtimes\langle t\rangle\bigr)
*_{H\times\langle t\rangle}
\bigl(K(S_2)\rtimes\langle t\rangle\bigr),
\]
which proves
\[
\mathcal G(G,H;S)
\cong
\mathcal G(G,H;S_1)
*_{H\times\mathbb Z}
\mathcal G(G,H;S_2).
\]

Now suppose that \(S\cong\mathbb Z/m\mathbb Z\). Then \(t^m\) fixes every factor \(G_s\) pointwise, and therefore
\[
\left(
\underset{s\in\mathbb Z/m\mathbb Z}{\ast_H}G
\right)\rtimes m\mathbb Z
\cong
\left(
\underset{s\in\mathbb Z/m\mathbb Z}{\ast_H}G
\right)\times\mathbb Z.
\]
This subgroup has index \(m\) in \(\mathcal G(G,H;S)\).

Finally, suppose that \(S\cong\mathbb Z\). Choosing the factor indexed by \(0\), all other factors are its conjugates by powers of \(t\). Thus
\[
\mathcal G(G,H;\mathbb Z)
=
\langle G,t\mid tht^{-1}=h,\ h\in H\rangle
\cong
G*_H(H\times\mathbb Z).
\]
\end{proof}

We now use Knopf's combination theorem for groups acting acylindrically on trees \cite{Knopf}. We first fix the convention on acylindricity.

\begin{definition}
Let \(k\geq 0\). An action of a group \(G\) on a simplicial tree \(T\) is called \emph{\(k\)-acylindrical} if the pointwise stabilizer of every geodesic segment of length \(k\) is finite. We call the action \emph{acylindrical} if it is \(k\)-acylindrical for some \(k\).
\end{definition}

This is the convention used by Delzant \cite[Definition, p. 1215]{Delzant} and Knopf \cite[Definition 3.1]{Knopf}; it differs slightly from Sela's original indexing convention \cite{sela} and from Bowditch's metric formulation \cite{Bowditch08}.

We shall use the following consequence of Knopf's theorem.

\begin{theorem}[Knopf {\cite[Corollary~4.4]{Knopf}}]
\label{thm:Knopf}
Let $G$ act acylindrically on a simplicial tree $T$. If every vertex stabilizer of the action $G\curvearrowright T$ belongs to $\mathcal C_{\mathrm{FJ}}$, then \(G\in\mathcal C_{\mathrm{FJ}}\).
\end{theorem}

Theorem~\ref{thm:Knopf} provides the main geometric tool for establishing FJ-admissibility. The difficulty is that the $\mathbb Z$--set occurring in Definition~\ref{def:FJ-admissible-pair} need not be transitive. For a single orbit the relevant groups admit simple Bass--Serre descriptions, whereas several orbits introduce additional amalgamations over \(H\times\mathbb Z\). We return to this issue in Subsubsection~\ref{subsubsec:concrete-FJ-pairs}.

\begin{theorem}[Local-to-global Farrell--Jones theorem]
\label{thm:FJ-graph-pairs}
Let $\Gamma$ be a finite simplicial graph and let \(\{(G_v,H_v)\}_{v\in V(\Gamma)}\) be a family of FJ-admissible pairs. Then \(\Pi_\Gamma(G_v,H_v)\in \mathcal C_{\mathrm{FJ}}\). In particular, the graph product of pairs satisfies the
Farrell--Jones conjecture.
\end{theorem}

\begin{proof}
We argue by induction on $|V(\Gamma)|$. If $\Gamma$ consists of a single vertex $v$, then \(\Pi_\Gamma(G_v,H_v)\cong G_v\), and $G_v\in\mathcal C_{\mathrm{FJ}}$ by Definition~\ref{def:FJ-admissible-pair}. Suppose that the result holds for graphs with fewer vertices. Choose $v\in V(\Gamma)$ and put
\[
P=\Pi_\Gamma(G_i,H_i),
\qquad
P'=\Pi_{\Gamma-\{v\}}(G_i,H_i).
\]
By the induction hypothesis, \(P'\in\mathcal C_{\mathrm{FJ}}\). Consider the pair-deletion epimorphism \(\rho_v\colon P\longrightarrow P'\). We first consider its kernel. By Proposition~\ref{prop:pair-deletion-kernel}, \(K_v:=\ker(\rho_v)\cong\underset{P'/R_v}{\ast_{H_v}}G_v\). Endow the nonempty set \(P'/R_v\) with the trivial \(\mathbb Z\)-action. Since \((G_v,H_v)\) is FJ-admissible, \(\mathcal G(G_v,H_v;P'/R_v)\cong K_v\times\mathbb Z\) belongs to \(\mathcal C_{\mathrm{FJ}}\). By subgroup closure, \(K_v\in\mathcal C_{\mathrm{FJ}}\). Now let \(C\leq P'\) be an infinite cyclic subgroup. By Proposition~\ref{prop:pair-deletion-cyclic}, \(\rho_v^{-1}(C)\cong\left(\underset{P'/R_v}{\ast_{H_v}}G_v\right)\rtimes C\). Since \(C\cong\mathbb Z\), its action permutes the factors indexed by \(P'/R_v\) and fixes the common subgroup \(H_v\) pointwise. Hence this is precisely one of the groups occurring in the definition of FJ-admissibility of \((G_v,H_v)\). Therefore \(\rho_v^{-1}(C)\in\mathcal C_{\mathrm{FJ}}\). The Farrell--Jones inheritance principle now yields \(P\in\mathcal C_{\mathrm{FJ}}\), as required.
\end{proof}

\begin{corollary}
\label{cor:FJ-polyhedral-pairs}
Let \(Z=Z_L(\underline A,\underline B)\) be a polyhedral product such that the associated actions \(G_i\curvearrowright E_i\) are transitive. Choose stabilizers $H_i\leq G_i$, so that \(E_i\cong G_i/H_i\). If every pair $(G_i,H_i)$ is FJ-admissible, then \(\pi_1(Z)\cong\Pi_{L_1}(G_i,H_i)\) satisfies the Farrell--Jones conjecture.
\end{corollary}

\begin{proof}
The fundamental-group description gives \(\pi_1(Z)\cong\Pi_{L_1}(G_i,H_i)\). The conclusion follows from Theorem~\ref{thm:FJ-graph-pairs}.
\end{proof}

\begin{remark}
No asphericity hypothesis is required in Corollary~\ref{cor:FJ-polyhedral-pairs}.  Asphericity is needed only when one passes from the group-theoretic Farrell--Jones statement to rigidity of the polyhedral product itself.
\end{remark}

\subsubsection{Concrete classes and alternative criteria}
\label{subsubsec:concrete-FJ-pairs}

Theorem~\ref{thm:FJ-graph-pairs} reduces the Farrell--Jones conjecture for graph products of pairs to the vertexwise FJ-admissibility condition.  We now discuss several situations in which the required Farrell--Jones conclusion follows from existing results. These cases arise from different mechanisms: the Gandini--R\"uping theorem for ordinary graph products, combination theorems for acylindrical actions on trees, extension arguments for normal pairs, and the Coxeter-group description of finite group pairs.

\paragraph{Ordinary graph products.}

The first case is the ordinary graph-product situation. If \(H=1\), then \(\underset{s\in S}{\ast_H}G=\ast_{s\in S}G\), and hence \(\mathcal G(G,1;S)=\left(\ast_{s\in S}G\right)\rtimes\mathbb Z\). Here the generator of \(\mathbb Z\) acts by permuting the free factors according to its action on \(S\).

Gandini--R\"uping prove precisely that this group satisfies the Farrell--Jones conjecture whenever \(G\) does \cite[Lemma~2.8]{GandiniRueping}. We therefore obtain the following.

\begin{prop}
\label{prop:trivial-pairs-FJ-admissible}
If \(G\in\mathcal C_{\mathrm{FJ}}\), then the pair \((G,1)\) is FJ-admissible.
\end{prop}

Consequently, Theorem~\ref{thm:FJ-graph-pairs} contains the Gandini--R\"uping theorem for ordinary graph products as the special case in which all \(H_i\) are trivial.

\medskip

\paragraph{Acylindrical pairs.}
The next source of FJ-admissible pairs comes from Bass--Serre theory. Let $H\leq G$ with \(G\in\mathcal C_{\mathrm{FJ}}\). For a transitive $\mathbb Z$--set $\mathcal O$, there are two possibilities. If \(\mathcal O\cong\mathbb Z/m\mathbb Z\), then \(\mathcal G(G,H;\mathcal O)\) contains \(\left(\underset{m}{\ast_H}G \right)\times\mathbb Z\) as a subgroup of index $m$. If \(\mathcal O\cong\mathbb Z\), then Lemma~\ref{lem:relative-orbit-decomposition} gives \(\mathcal G(G,H;\mathbb Z) \cong G*_H(H\times\mathbb Z)\). This leads to the following criterion.

\begin{prop}[Single-orbit acylindrical criterion]
\label{prop:single-orbit-FJ}
Let \(H\leq G\) and suppose that \(G\in\mathcal C_{\mathrm{FJ}}\). Let \(\mathcal O\) be a transitive \(\mathbb Z\)-set.

\begin{enumerate}
\item
If \(\mathcal O\cong\mathbb Z/m\mathbb Z\) and the natural Bass--Serre action of \(\underset{m}{\ast_H}G\) is acylindrical, then \(\mathcal G(G,H;\mathcal O)\in\mathcal C_{\mathrm{FJ}}\).

\item
If \(\mathcal O\cong\mathbb Z\) and the natural Bass--Serre action of \(G*_H(H\times\mathbb Z)\) is acylindrical, then \(\mathcal G(G,H;\mathcal O) \in\mathcal C_{\mathrm{FJ}}\).
\end{enumerate}
\end{prop}

\begin{proof}
For~(1), the Bass--Serre tree of \(\underset{m}{\ast_H}G\) has vertex stabilizers conjugate to copies of \(G\) and to \(H\). Both belong to \(\mathcal C_{\mathrm{FJ}}\): this is an assumption for \(G\), and follows for \(H\) from subgroup closure. Knopf's theorem therefore gives \(\underset{m}{\ast_H}G \in\mathcal C_{\mathrm{FJ}}\). Closure of \(\mathcal C_{\mathrm{FJ}}\) under finite direct products then gives \(\left(
\underset{m}{\ast_H}G
\right)\times\mathbb Z
\in\mathcal C_{\mathrm{FJ}}\). Since this group has finite index in \(\mathcal G(G,H;\mathcal O)\), closure under finite-index overgroups proves the assertion.

For~(2), since \(H\leq G\), subgroup closure gives \(H\in\mathcal C_{\mathrm{FJ}}\), and hence \(H\times\mathbb Z\in\mathcal C_{\mathrm{FJ}}\). By Lemma~\ref{lem:relative-orbit-decomposition}, \(\mathcal G(G,H;\mathbb Z) \cong G*_H(H\times\mathbb Z)\). The result follows from Knopf's theorem applied to the assumed acylindrical Bass--Serre action.
\end{proof}

For several orbits an additional difficulty appears. If \(S=S_1\sqcup S_2\) is a decomposition into invariant subsets, then Lemma~\ref{lem:relative-orbit-decomposition} gives \(\mathcal G(G,H;S)\cong\mathcal G(G,H;S_1)*_{H\times\mathbb Z}\mathcal G(G,H;S_2)\). Thus the passage from one orbit to an arbitrary $\mathbb Z$--set introduces amalgamations over $H\times\mathbb Z$. Their Bass--Serre actions need not be acylindrical automatically.

One convenient sufficient condition is the following.

\begin{prop}[Acylindrical FJ-admissibility criterion]
\label{prop:acylindrical-FJ-admissible}
Let \(H\leq G\), with \(G\in\mathcal C_{\mathrm{FJ}}\). Assume that:
\begin{enumerate}
\item for every transitive \(\mathbb Z\)-set \(\mathcal O\), the group \(\mathcal G(G,H;\mathcal O)\) belongs to \(\mathcal C_{\mathrm{FJ}}\);
\item whenever \(S_1\) and \(S_2\) are nonempty finite unions of \(\mathbb Z\)-orbits, the Bass--Serre action associated with \(\mathcal G(G,H;S_1)*_{H\times\mathbb Z}\mathcal G(G,H;S_2)\) is acylindrical.
\end{enumerate}
Then \((G,H)\) is FJ-admissible.
\end{prop}

\begin{proof}
By Lemma~\ref{lem:relative-orbit-decomposition} and Theorem~\ref{thm:Knopf}, induction on the number of orbits shows that \(\mathcal G(G,H;S)\in\mathcal C_{\mathrm{FJ}}\) whenever \(S\) is a finite union of \(\mathbb Z\)-orbits.

For an arbitrary \(\mathbb Z\)-set \(S\), let \(\mathcal S\) denote the directed set of finite unions of orbits in \(S\). Then \(\mathcal G(G,H;S)\cong\operatorname*{colim}_{S_0\in\mathcal S}\mathcal G(G,H;S_0)\). Since \(\mathcal C_{\mathrm{FJ}}\) is closed under directed colimits, \(\mathcal G(G,H;S)\in\mathcal C_{\mathrm{FJ}}\). Hence \((G,H)\) is FJ-admissible.
\end{proof}

\begin{remark}
The preceding proposition is deliberately formulated as a criterion rather than as a classification of FJ-admissible pairs. In concrete applications one may verify acylindricity through almost malnormality or other geometric hypotheses on the inclusion \(H\leq G\). Recall that a subgroup \(H\leq G\) is \emph{almost malnormal} if \(H\cap gHg^{-1}\) is finite for every \(g\in G\setminus H\) (see \cite[Definition 3.3]{Knopf}). Determining the largest natural class of pairs for which these conditions hold remains an interesting problem.
\end{remark}

\paragraph{Free-factor pairs.}

There is also a non-normal class which can be treated directly using the free-product combination results of Gandini--R\"uping \cite{GandiniRueping}.

\begin{prop}[Free-factor pairs]
\label{prop:free-factor-FJ-admissible}
Suppose that \(G=H*K\) and that \(H,K\in\mathcal C_{\mathrm{FJ}}\). Then the pair \((G,H)\) is FJ-admissible.
\end{prop}
\begin{proof}
Since \(G=H*K\) and \(H,K\in\mathcal C_{\mathrm{FJ}}\), closure under free products gives \(G\in\mathcal C_{\mathrm{FJ}}\). Let \(S\) be a nonempty \(\mathbb Z\)-set. Since \(G=H*K\), amalgamating copies of \(G\) along their common factor \(H\) gives \(\underset{s\in S}{\ast_H}G\cong H*\left(\underset{s\in S}{\ast}K\right)\). Under this identification, the \(\mathbb Z\)-action fixes \(H\) pointwise and permutes the copies of \(K\) according to its action on \(S\). Put \(K_S=\underset{s\in S}{\ast}K\). By \cite[Lemma~2.8]{GandiniRueping}, \(K_S\rtimes\mathbb Z \in\mathcal C_{\mathrm{FJ}}\). Moreover, \(H\rtimes_{\mathrm{id}}\mathbb Z=H\times\mathbb Z\in\mathcal C_{\mathrm{FJ}}\). Applying \cite[Lemma~2.7]{GandiniRueping} to the free product \(H*K_S\) and the automorphism which is the identity on \(H\) and the given permutation automorphism on \(K_S\), we obtain \(\left(H*K_S\right)\rtimes\mathbb Z\in\mathcal C_{\mathrm{FJ}}\). Hence \(\mathcal G(G,H;S)\in\mathcal C_{\mathrm{FJ}}\). Since \(S\) was arbitrary, \((G,H)\) is FJ-admissible.
\end{proof}

\begin{example}
\label{ex:F2-cyclic-pair}
Let \(G=F(a,b)=\langle a\rangle *\langle b\rangle\) and \(H=\langle a\rangle\). Since both \(\langle a\rangle\cong\mathbb Z\) and \(\langle b\rangle\cong\mathbb Z\) belong to \(\mathcal C_{\mathrm{FJ}}\), Proposition~\ref{prop:free-factor-FJ-admissible} shows that \((F_2,\langle a\rangle)\) is FJ-admissible. In particular, this gives an infinite non-normal example of an FJ-admissible pair.
\end{example}

\paragraph{Normal pairs and cyclic extensions.}

The normal case admits a stronger formulation which isolates the precise local property needed in the pair-deletion argument. Namely, what is required is not virtual solvability, or even amenability itself, but stability of the Farrell--Jones conjecture under cyclic extensions of the common subgroup.

\begin{definition}
\label{def:FJ-cyclic-extension-stable}
A group \(H\) will be called \emph{FJ-stable under cyclic extensions} if \(H\in\mathcal C_{\mathrm{FJ}}\) and, for every automorphism \(\alpha\in\operatorname{Aut}(H)\), the semidirect product \(H\rtimes_\alpha\mathbb Z\) belongs to \(\mathcal C_{\mathrm{FJ}}\).
\end{definition}

\begin{remark}
The condition in
Definition~\ref{def:FJ-cyclic-extension-stable} is tailored to the Farrell--Jones inheritance principle.  Indeed, whenever
\[
1\longrightarrow H\longrightarrow E\longrightarrow\mathbb Z
\longrightarrow1
\]
is exact, the extension splits and \(E\cong H\rtimes_\alpha\mathbb Z\) for some \(\alpha\in\operatorname{Aut}(H)\). Thus FJ-cyclic-extension stability is precisely the condition required to control all infinite cyclic preimages of a quotient whose kernel is \(H\).
\end{remark}

\begin{remark}
The appearance of cyclic extensions is not accidental. Wu's proof of the Farrell--Jones conjecture for elementary amenable groups of finite Hirsch length likewise reduces, through the inheritance principle, to groups which are extensions by an infinite cyclic group. This is precisely the type of subgroup produced by the cyclic preimages in the pair-deletion argument.
\end{remark}

\begin{prop}[Normal FJ-admissible pairs]
\label{prop:normal-FJ-admissible}
Let \(H\trianglelefteq G\) and put \(Q=G/H\). Assume that \(Q\in\mathcal C_{\mathrm{FJ}}\) and that \(H\) is FJ-stable under cyclic extensions. Then the pair \((G,H)\) is FJ-admissible.
\end{prop}
\begin{proof}
We first show that \(G\in\mathcal C_{\mathrm{FJ}}\). Consider the quotient homomorphism \(q\colon G\longrightarrow Q\). By assumption, \(Q\in\mathcal C_{\mathrm{FJ}}\). The kernel of \(q\) is \(H\), which belongs to \(\mathcal C_{\mathrm{FJ}}\) by the assumption that \(H\) is FJ-stable under cyclic extensions. Let \(C\leq Q\) be an infinite cyclic subgroup. Then
\[
1\longrightarrow H
\longrightarrow q^{-1}(C)
\longrightarrow C
\longrightarrow1
\]
is exact. Since \(C\cong\mathbb Z\) is free, this extension splits, and hence \(q^{-1}(C)\cong H\rtimes_\alpha\mathbb Z\) for some \(\alpha\in\operatorname{Aut}(H)\). Therefore \(q^{-1}(C)\in\mathcal C_{\mathrm{FJ}}\). The Farrell--Jones inheritance principle now gives \(G\in\mathcal C_{\mathrm{FJ}}\).

Let \(S\) be a nonempty \(\mathbb Z\)--set. Since \(H\trianglelefteq G\), the common subgroup \(H\) is normal in \(\underset{s\in S}{\ast_H}G\). Moreover, the \(\mathbb Z\)--action permutes the factors and fixes \(H\) pointwise. Hence there is a short exact sequence
\[
1\longrightarrow H
\longrightarrow
\mathcal G(G,H;S)
\stackrel{\overline q}{\longrightarrow}
\left(
\underset{s\in S}{\ast}Q
\right)\rtimes\mathbb Z
\longrightarrow1.
\]

By \cite[Lemma~2.8]{GandiniRueping}, \(\left(\underset{s\in S}{\ast}Q \right)\rtimes\mathbb Z\in\mathcal C_{\mathrm{FJ}}\). Again, \(\ker(\overline q)=H\in\mathcal C_{\mathrm{FJ}}\). Let \(C\) be an infinite cyclic subgroup of the quotient. Then
\[
1\longrightarrow H
\longrightarrow
\overline q^{-1}(C)
\longrightarrow
C
\longrightarrow1
\]
is exact, and hence \(\overline q^{-1}(C) \cong H\rtimes_\beta\mathbb Z\) for some automorphism \(\beta\in\operatorname{Aut}(H)\). By the assumption that \(H\) is FJ-stable under cyclic extensions, \(\overline q^{-1}(C)\in\mathcal C_{\mathrm{FJ}}\). The inheritance principle therefore gives \(\mathcal G(G,H;S)\in\mathcal C_{\mathrm{FJ}}\). Since this holds for every nonempty \(\mathbb Z\)--set \(S\), the pair \((G,H)\) is FJ-admissible.
\end{proof}

\begin{prop}[Normal pairs with Farrell--Jones quotient]
\label{prop:normal-pairs-G-Q-FJ}
Let \(H\trianglelefteq G\) and put \(Q=G/H\). Assume that \(G,Q\in\mathcal C_{\mathrm{FJ}}\). Then the pair \((G,H)\) is FJ-admissible.
\end{prop}
\begin{proof}
Let \(S\) be a nonempty \(\mathbb Z\)-set. Since \(H\trianglelefteq G\), there is a short exact sequence
\[
1\longrightarrow H
\longrightarrow
\mathcal G(G,H;S)
\stackrel{\overline q}{\longrightarrow}
\left(
\underset{s\in S}{\ast}Q
\right)\rtimes\mathbb Z
\longrightarrow1.
\]

By \cite[Lemma~2.8]{GandiniRueping}, \(\left(\underset{s\in S}{\ast}Q\right)\rtimes\mathbb Z\in\mathcal C_{\mathrm{FJ}}\). Moreover, \(H\in\mathcal C_{\mathrm{FJ}}\) by subgroup closure, since \(H\leq G\). Let \(C\) be an infinite cyclic subgroup of the quotient, and let \(c\) be a generator of \(C\). Choose a lift \(t=k\tau^n\in\mathcal G(G,H;S)\), where \(k\in\underset{s\in S}{\ast_H}G\) and \(\tau\) denotes the generator of the acting copy of \(\mathbb Z\). Let \(\mathrm{fold}\colon\underset{s\in S}{\ast_H}G\longrightarrow G\) be the homomorphism which is the identity on every copy of \(G\), and put \(g=\mathrm{fold}(k)\). Since \(\tau\) fixes \(H\) pointwise, conjugation by \(t\) on \(H\) coincides with conjugation by \(g\). Consequently, \(\overline q^{-1}(C) \cong H\rtimes_{c_g}\mathbb Z\). The homomorphism
\[
H\rtimes_{c_g}\mathbb Z
\longrightarrow
G\times\mathbb Z,
\qquad
(h,m)\longmapsto (hg^m,m),
\]
is injective. Hence \(\overline q^{-1}(C)\leq G\times\mathbb Z\). Since \(G\times\mathbb Z\in\mathcal C_{\mathrm{FJ}}\), subgroup closure gives \(\overline q^{-1}(C)\in\mathcal C_{\mathrm{FJ}}\). The Farrell--Jones inheritance principle therefore yields \(\mathcal G(G,H;S)\in\mathcal C_{\mathrm{FJ}}\). Since this holds for every nonempty \(\mathbb Z\)-set \(S\), the pair \((G,H)\) is FJ-admissible.
\end{proof}

\begin{remark}
Propositions~\ref{prop:normal-FJ-admissible} and \ref{prop:normal-pairs-G-Q-FJ} are complementary. The former does not assume \(G\in\mathcal C_{\mathrm{FJ}}\), but instead derives this from the cyclic-extension stability of \(H\). The latter requires \(G\in\mathcal C_{\mathrm{FJ}}\) from the outset, but imposes no additional condition on \(H\).
\end{remark}

\begin{corollary}[Normal virtually solvable pairs]
\label{cor:normal-virtually-solvable-admissible}
Let \(H\trianglelefteq G\), \(Q=G/H\). Suppose that \(H\) is virtually solvable and \(Q\in\mathcal C_{\mathrm{FJ}}\). Then \((G,H)\) is FJ-admissible.
\end{corollary}
\begin{proof}
By Wegner's theorem \cite{Wegner}, \(H\in\mathcal C_{\mathrm{FJ}}\). Moreover, for every \(\alpha\in\operatorname{Aut}(H)\), the semidirect product \(H\rtimes_\alpha\mathbb Z\) is again virtually solvable, and hence belongs to \(\mathcal C_{\mathrm{FJ}}\) by Wegner's theorem. Thus \(H\) is FJ-stable under cyclic extensions. The conclusion follows from Proposition~\ref{prop:normal-FJ-admissible}.
\end{proof}

\begin{corollary}[Normal elementary amenable pairs]
\label{cor:normal-elementary-amenable-admissible}
Let \(H\trianglelefteq G\), \(Q=G/H\). Suppose that \(H\) is elementary amenable of finite Hirsch length and \(Q\in\mathcal C_{\mathrm{FJ}}\). Then \((G,H)\) is FJ-admissible.
\end{corollary}

\begin{proof}
By \cite[Theorem~A]{WuElementaryAmenableFJC}, the group \(H\) satisfies the Farrell--Jones conjecture with finite wreath products. For every automorphism \(\alpha\in\operatorname{Aut}(H)\), the semidirect product \(H\rtimes_\alpha\mathbb Z\) is again elementary amenable. Moreover, \(h(H\rtimes_\alpha\mathbb Z)\leq h(H)+1<\infty\). Hence \(H\rtimes_\alpha\mathbb Z\in\mathcal C_{\mathrm{FJ}}\) again by \cite[Theorem~A]{WuElementaryAmenableFJC}. Thus \(H\) is FJ-stable under cyclic extensions, and the conclusion follows from Proposition~\ref{prop:normal-FJ-admissible}.
\end{proof}

\begin{remark}[Houghton-type examples]
\label{rem:Houghton-pairs}
Corollary~\ref{cor:normal-elementary-amenable-admissible} is genuinely stronger than the corresponding virtually solvable statement. Wu \cite{WuElementaryAmenableFJC} points out that Houghton groups give elementary amenable groups of finite Hirsch length which are not virtually solvable. For example, let \(H\) be such a Houghton group and let \(Q\in\mathcal C_{\mathrm{FJ}}\). Then, for \(G=H\times Q\), the inclusion \(H\trianglelefteq G\) has quotient \(G/H\cong Q\). Hence \((G,H)\) is FJ-admissible by Corollary~\ref{cor:normal-elementary-amenable-admissible}. Thus Theorem~\ref{thm:FJ-graph-pairs} produces graph products of pairs satisfying the Farrell--Jones conjecture whose local normal subgroups need not be virtually solvable.
\end{remark}

\medskip

\noindent\emph{An alternative global normal-pair criterion.}
When all subgroups \(H_i\) are normal, one can also work directly with the global extension, without verifying FJ-admissibility vertex by vertex.

Suppose that \(H_i\trianglelefteq G_i\) for every \(i\), and put \(Q_i=G_i/H_i\). By Proposition~\ref{prop:normal-pair-extension}, there is a short exact sequence
\[
1
\longrightarrow
N
\longrightarrow
\Pi_{L_1}(G_i,H_i)
\stackrel{\rho}{\longrightarrow}
Q
\longrightarrow
1,
\]
where
\[
N=\prod_i H_i,
\qquad
Q=\mathop{\ast}_{L_1}Q_i.
\]

\begin{prop}[Global normal-pair criterion]
\label{prop:FJ-normal-pairs}
Assume that:
\begin{enumerate}
\item every \(Q_i\) belongs to \(\mathcal C_{\mathrm{FJ}}\);
\item \(N=\prod_iH_i\) belongs to \(\mathcal C_{\mathrm{FJ}}\);
\item for every infinite cyclic subgroup \(C\leq Q\), \(\rho^{-1}(C)\in\mathcal C_{\mathrm{FJ}}\).
\end{enumerate}
Then \(\Pi_{L_1}(G_i,H_i)\in\mathcal C_{\mathrm{FJ}}\).
\end{prop}
\begin{proof}
By Theorem~\ref{thm:GR-graph-products}, \(Q=\mathop{\ast}_{L_1}Q_i \in\mathcal C_{\mathrm{FJ}}\). The kernel of \(\rho\) is \(N\), which belongs to \(\mathcal C_{\mathrm{FJ}}\) by assumption. By assumption, the preimage \(\rho^{-1}(C)\) also belongs to \(\mathcal C_{\mathrm{FJ}}\) for every infinite cyclic subgroup \(C\leq Q\). The conclusion therefore follows from the Farrell--Jones inheritance principle.
\end{proof}

\begin{remark}[Amenability and the Farrell--Jones problem.]
The role of amenability in the present section should be distinguished from that of acylindricity. The latter provides a geometric combination mechanism through Knopf's theorem, whereas amenability identifies natural classes of local subgroups which may already satisfy the Farrell--Jones conjecture. This distinction is particularly relevant in view of the amenable-category results obtained earlier in the paper. The same local subgroups \(H_i\leq G_i\) therefore enter the theory in two logically distinct ways: amenability controls amenable category and simplicial-volume phenomena, while Farrell--Jones and Bass--Serre properties of the inclusions control assembly maps and rigidity. The Farrell--Jones conjecture is not known for arbitrary amenable groups. However, virtually solvable groups are covered by Wegner's theorem, and elementary amenable groups of finite Hirsch length satisfy FJCw by \cite[Theorem~A]{WuElementaryAmenableFJC}.
\end{remark}

Summarizing the preceding results and Theorem~\ref{thm:FJ-graph-pairs}, we obtain the following.

\begin{theorem}[Farrell--Jones combination theorem for graph products of pairs] \label{thm:main-FJ-graph-pairs}
Let $\Gamma$ be a finite simplicial graph and let \(\{(G_v,H_v)\}_{v\in V(\Gamma)}\) be a family of group pairs. Assume that, for every \(v\in V(\Gamma)\), at least one of the following conditions holds:
\begin{enumerate}[(i)]
\item \(H_v=\{1\}\) and \(G_v\in\mathcal C_{\mathrm{FJ}}\);
\item \(H_v\trianglelefteq G_v\), \(G_v/H_v\in\mathcal C_{\mathrm{FJ}}\), and either \(G_v\in\mathcal C_{\mathrm{FJ}}\) or \(H_v\) is FJ-stable under cyclic extensions;
\item \(G_v=H_v*K_v\) with \(H_v,K_v\in\mathcal C_{\mathrm{FJ}}\);
\item \(G_v\in\mathcal C_{\mathrm{FJ}}\) and the pair \((G_v,H_v)\) satisfies the hypotheses of Proposition~\ref{prop:acylindrical-FJ-admissible}.
\end{enumerate}

Then \(\Pi_\Gamma(G_v,H_v)\in\mathcal C_{\mathrm{FJ}}\).
\end{theorem}

\begin{remark}
\label{rem:concrete-normal-main-theorem}
Condition~{\rm(ii)} of Theorem~\ref{thm:main-FJ-graph-pairs} contains two important concrete classes. If \(H_v\) is virtually solvable, then \(H_v\) is FJ-stable under cyclic extensions by Wegner's theorem. If \(H_v\) is elementary amenable of finite Hirsch length, then \(H_v\) is FJ-stable under cyclic extensions by \cite[Theorem~A]{WuElementaryAmenableFJC}. Thus, provided \(G_v/H_v\in\mathcal C_{\mathrm{FJ}}\), both classes are covered by condition~{\rm(ii)}.
\end{remark}

The importance of Theorem~\ref{thm:main-FJ-graph-pairs} is that its hypotheses are vertexwise: different vertices may be treated by different mechanisms. Thus a single graph product of pairs may simultaneously contain ordinary graph-product vertices, free-factor pairs, normal pairs treated by cyclic-extension methods, and non-normal pairs controlled by acylindrical Bass--Serre actions.

\begin{corollary}[Farrell--Jones for polyhedral products]
\label{cor:main-FJ-polyhedral-products}
Let \(Z=Z_L(\underline A,\underline B)\) be a polyhedral product such that each action \(G_i\curvearrowright E_i\) is transitive. Choose \(H_i\leq G_i\) such that \(E_i\cong G_i/H_i\). Suppose that every pair \((G_i,H_i)\) satisfies one of the conditions of Theorem~\ref{thm:main-FJ-graph-pairs}. Then \(\pi_1(Z)\cong\Pi_{L_1}(G_i,H_i)\in\mathcal C_{\mathrm{FJ}}\).
\end{corollary}

\medskip

Theorem~\ref{thm:main-FJ-graph-pairs} collects the classes for which the local FJ-admissibility condition can be verified directly. There is also an independent global class which does not require this local criterion.

\begin{prop}
\label{prop:FJ-finite-pairs}
Let \(G=\Pi_\Gamma(G_i,H_i)\) be a graph product of pairs in which every \(G_i\) is finite. Then \(G\in\mathcal C_{\mathrm{FJ}}\).
\end{prop}

\begin{proof}
Januszkiewicz--Świ\k{a}tkowski \cite[Corollary~5.15]{JanuszkiewiczSwiatkowski} show that graph products of finite group pairs are commensurable with Coxeter groups. Coxeter groups are \(\operatorname{CAT}(0)\)-groups and hence belong to \(\mathcal C_{\mathrm{FJ}}\); see~\cite[Proposition~1.1]{GandiniRueping}.
\end{proof}

\begin{remark}
The finite-pair case illustrates that FJ-admissibility is a sufficient local condition for the pair-deletion theorem, but is not intended to be a necessary condition for a graph product of pairs to satisfy Farrell--Jones. In some classes, such as finite group pairs, the global group can be handled by an independent geometric or commensurability argument.
\end{remark}

\medskip

The preceding examples show that the Farrell--Jones problem for graph products of pairs admits several complementary approaches. The pair-deletion method is local and applies without normality, provided the vertex pairs are FJ-admissible. Free-factor decompositions give a direct non-normal class through the combination results of Gandini--R\"uping, while acylindrical Bass--Serre actions provide a geometric mechanism for treating more general inclusions. For normal pairs, cyclic-extension and quotient arguments give a different route, whereas finite group pairs can be handled globally through their relation with Coxeter groups.

We close this subsection with the following problem. 

\begin{question}
\label{ques:FJ-admissible-pairs} 
Which natural classes of inclusions \(H\leq G\) are FJ-admissible beyond the normal, free-factor, and acylindrical classes considered above?

More specifically, can the acylindricity hypotheses in Proposition~\ref{prop:acylindrical-FJ-admissible} be replaced by an appropriate relative finite-homotopy \(\mathcal F\)-amenability condition for the Bass--Serre actions arising from pair deletion?
\end{question}

The preceding results reveal a local-to-global mechanism for rigidity of polyhedral products.  The local data are the vertex pairs \((G_i,H_i)\), which arise from the transitive actions \(G_i\curvearrowright \pi_0\bigl(p_i^{-1}(B_i)\bigr)\). Pair deletion converts the global graph product of pairs into a sequence of amalgamated free-product problems whose cyclic preimages are controlled by the corresponding vertex pair.

Thus the Farrell--Jones problem reduces to determining natural classes of FJ-admissible pairs.  Any such class immediately yields, through Theorem~\ref{thm:FJ-graph-pairs}, a Farrell--Jones theorem for the associated graph products of pairs, and hence for the fundamental groups of the corresponding polyhedral products.

For closed aspherical manifold polyhedral products of dimension at least $5$, Theorem~\ref{thm:Borel-FJ-polyhedral-products} then converts these group-theoretic results into topological rigidity statements.

\subsection{The Borel conjecture}
\label{subsec:Borel}

Recall that a closed manifold \(M\) is called \emph{topologically rigid} if every homotopy equivalence \(f\colon N\longrightarrow M\) from a closed manifold \(N\) is homotopic to a homeomorphism. The Borel conjecture asserts that every closed aspherical manifold is topologically rigid.

The connection with the Farrell--Jones conjecture is provided by surgery theory. We shall use the following standard consequence of the work of Bartels and L\"uck~\cite{BartelsLuck}.

\begin{theorem}[Bartels--L\"uck {\cite{BartelsLuck}}]
\label{thm:BL-Borel}
Let \(G\) be torsion-free. Suppose that the \(K\)-theoretic Farrell--Jones assembly map is an isomorphism in degrees \(m\leq 0\) and surjective in degree \(m=1\), and that the \(L\)-theoretic Farrell--Jones assembly map is an isomorphism in all degrees, including all orientation twists. Then the Borel conjecture holds for closed aspherical manifolds of dimension at least \(5\) with fundamental group \(G\).
\end{theorem}

In particular, the form of the Farrell--Jones conjecture encoded by \(\mathcal C_{\mathrm{FJ}}\) is sufficient for the preceding theorem.

\begin{corollary}[Borel rigidity for polyhedral products]
\label{thm:Borel-FJ-polyhedral-products}
Let \(M=Z_L(\underline A,\underline B)\) be a closed aspherical manifold polyhedral product of dimension at least \(5\). If \(\pi_1(M)\in\mathcal C_{\mathrm{FJ}}\), then \(M\) satisfies the Borel conjecture.
\end{corollary}

\begin{proof}
Since \(M\) is a closed aspherical manifold, \(\pi_1(M)\) is torsion-free by Remark~\ref{rem:torsion-free-manifold-group}. The Farrell--Jones hypothesis implies the assembly-map assumptions of Theorem~\ref{thm:BL-Borel}. Hence \(M\) is topologically rigid.
\end{proof}

Combining this with the Farrell--Jones combination theorem for graph products of pairs gives the following direct application.

\begin{corollary}[Borel rigidity from graph products of pairs]
\label{cor:Borel-main-graph-pairs}
Let \(M=Z_L(\underline A,\underline B)\) be a closed aspherical manifold polyhedral product of dimension at least \(5\). Suppose that each component action \(G_i\curvearrowright E_i\) is transitive, and choose \(H_i\leq G_i\) such that \(E_i\cong G_i/H_i\). Assume that, for every \(i\), the pair \((G_i,H_i)\) satisfies one of the conditions of Theorem~\ref{thm:main-FJ-graph-pairs}. Then \(M\) satisfies the Borel conjecture.
\end{corollary}

\begin{proof}
By Corollary~\ref{cor:main-FJ-polyhedral-products}, \(\pi_1(M) \cong \Pi_{L_1}(G_i,H_i) \in\mathcal C_{\mathrm{FJ}}\). The conclusion follows from Corollary~\ref{thm:Borel-FJ-polyhedral-products}.
\end{proof}

\begin{remark}
The preceding corollary applies vertexwise. Thus different local pairs may simultaneously be ordinary graph-product pairs, free-factor pairs, normal pairs covered by the cyclic-extension criteria, or pairs satisfying the acylindrical Bass--Serre hypotheses of Proposition~\ref{prop:acylindrical-FJ-admissible}.
\end{remark}

We close with several manifold examples  illustrating the scope of the preceding results. They show that the local pairs considered above arise naturally from manifold pairs and provide useful consistency checks for the Farrell--Jones and rigidity machinery. The real moment-angle case represents the basic non-transitive situation. M\"obius-type polyhedral products give a normal but non-split example, whereas surface-boundary pairs provide natural non-normal inclusions which are not covered by the elementary normal or free-factor criteria.

\begin{example}[Real moment-angle manifolds
{cf.~\cite{gkeneralis-real}}]
\label{ex:Borel-real-MAM}
Let \(L\) be a finite flag triangulation of a sphere and suppose that \(\mathcal Z_L^{\mathbb R}=Z_L(D^1,S^0)\) has dimension at least \(5\). Its fundamental group is a torsion-free finite-index subgroup of the right-angled Coxeter group \(W_L\). Since \(W_L\) is a \(\operatorname{CAT}(0)\) group, it belongs to \(\mathcal C_{\mathrm{FJ}}\), and hence so does \(\pi_1(\mathcal Z_L^{\mathbb R})\) by subgroup closure. Therefore \(\mathcal Z_L^{\mathbb R}\) satisfies the Borel conjecture.
\end{example}

Notice that the above example lies outside the transitive framework of Corollary~\ref{cor:Borel-main-graph-pairs}: here \(G_i=1\) acts on the two-point set \(E_i=S^0\), and the action is not transitive. Thus real moment-angle manifolds provide a basic model for the more general non-transitive relative graph-product situation.

\begin{example}[M\"obius-type polyhedral products]
Let \((A,B)=(\mathrm{M\ddot{o}},S^1)\), where \(\mathrm{M\ddot{o}}\) denotes the M\"obius band and \(S^1=\partial\mathrm{M\ddot{o}}\). Then \(\pi_1(\mathrm{M\ddot{o}})\cong\mathbb Z\), and \(\pi_1(S^1)\cong 2\mathbb Z\leq\mathbb Z\). In particular, the inclusion \(\partial\mathrm{M\ddot{o}}\hookrightarrow\mathrm{M\ddot{o}}\) induces
\[
\pi_1(\partial\mathrm{M\ddot{o}})\cong\mathbb Z
\longrightarrow
\pi_1(\mathrm{M\ddot{o}})\cong\mathbb Z
\]
with image \(2\mathbb Z\). Thus the associated local pair is \((\mathbb Z,2\mathbb Z)\). It is a normal but non-split pair: the extension
\[
1\longrightarrow 2\mathbb Z
\longrightarrow \mathbb Z
\longrightarrow \mathbb Z/2
\longrightarrow1
\]
does not split. Since \(2\mathbb Z\) is virtually solvable and \(\mathbb Z/2\in\mathcal C_{\mathrm{FJ}}\), Corollary~\ref{cor:normal-virtually-solvable-admissible} shows that \((\mathbb Z,2\mathbb Z)\) is FJ-admissible. Hence any closed aspherical manifold polyhedral product built from such local data and satisfying the hypotheses above has fundamental group in \(\mathcal C_{\mathrm{FJ}}\), and therefore satisfies the Borel conjecture in dimensions at least \(5\).
\end{example}

The above example should be viewed as a consistency check for the normal-pair machinery rather than as a new Farrell--Jones case, since the resulting groups are also accessible by finite-cover arguments.

\begin{example}[Surface-boundary pairs]
Let \(\Sigma_{g,1}\) be a compact orientable surface of genus \(g\) with one boundary component. Then
\[
G=\pi_1(\Sigma_{g,1})\cong F_{2g},
\qquad
H=\pi_1(\partial\Sigma_{g,1})\cong\mathbb Z.
\]
The subgroup \(H\leq G\) is generally neither normal nor a free factor. Thus this pair lies outside the elementary normal and free-factor criteria considered above. It therefore provides a natural test case for the acylindrical criterion and, more generally, for Question~\ref{ques:FJ-admissible-pairs}. Determining whether such boundary inclusions are FJ-admissible would produce further Farrell--Jones and rigidity results for polyhedral products built from manifolds with boundary.
\end{example}

\printbibliography

@article{davis_asphericity,
 author = {Davis, Michael W.},
 title = {Right-angularity, flag complexes, asphericity},
 fjournal = {Geometriae Dedicata},
 journal = {Geom. Dedicata},
 issn = {0046-5755},
 volume = {159},
 pages = {239--262},
 year = {2012},
 language = {English},
 doi = {10.1007/s10711-011-9654-4},
 zbMATH = {6062604},
 Zbl = {1293.57003}
}

@article{davis_kropholler_corrigenda,
 author = {Davis, Michael W. and Kropholler, Peter H.},
 title = {Criteria for asphericity of polyhedral products: corrigenda to ``{Right}-angularity, flag complexes, asphericity''},
 fjournal = {Geometriae Dedicata},
 journal = {Geom. Dedicata},
 issn = {0046-5755},
 volume = {179},
 pages = {39--44},
 year = {2015},
 language = {English},
 doi = {10.1007/s10711-015-0066-8},
 zbMATH = {6521447},
 Zbl = {1332.57003}
}

@book{davis_coxeter,
 author = {Davis, Michael W.},
 title = {The geometry and topology of {Coxeter} groups},
 edition = {2nd edition},
 fseries = {Springer Monographs in Mathematics},
 series = {Springer Monogr. Math.},
 issn = {1439-7382},
 %isbn = {978-3-031-91302-0; 978-3-031-91303-7},
 year = {2025},
 publisher = {Cham: Springer},
 language = {English},
 doi = {10.1007/978-3-031-91303-7},
 keywords = {20-02,57-02,20F55,20F65,57M07,51F15,20C08},
 zbMATH = {8014470},
 Zbl = {1573.20001}
}

@book{tt,
 author = {Buchstaber, Victor M. and Panov, Taras E.},
 title = {Toric topology},
 fseries = {Mathematical Surveys and Monographs},
 series = {Math. Surv. Monogr.},
 issn = {0076-5376},
 volume = {204},
 isbn = {978-1-4704-2214-1},
 year = {2015},
 publisher = {Providence, RI: American Mathematical Society (AMS)},
 language = {English},
 zbMATH = {6471875},
 Zbl = {1375.14001}
}

@article{bbcg2010,
 author = {Bahri, A. and Bendersky, M. and Cohen, F. R. and Gitler, S.},
 title = {The polyhedral product functor: {A} method of decomposition for moment-angle complexes, arrangements and related spaces},
 fjournal = {Advances in Mathematics},
 journal = {Adv. Math.},
 issn = {0001-8708},
 volume = {225},
 number = {3},
 pages = {1634--1668},
 year = {2010},
 language = {English},
 doi = {10.1016/j.aim.2010.03.026},
 zbMATH = {5784902},
 Zbl = {1197.13021}
}

@article{bbcg2009,
 author = {Bahri, A. and Bendersky, M. and Cohen, F. R. and Gitler, S.},
 title = {Gitler decompositions of the polyhedral product functor with applications to moment-angle complexes and related spaces},
 fjournal = {Proceedings of the National Academy of Sciences of the United States of America},
 journal = {Proc. Natl. Acad. Sci. USA},
 issn = {0027-8424},
 volume = {106},
 number = {30},
 pages = {12241--12244},
 year = {2009},
 language = {English},
 doi = {10.1073/pnas.0905159106},
 url = {europepmc.org/articles/pmc2713394},
 zbMATH = {5840882},
 Zbl = {1203.57009}
}

@misc{gken_rigidityCH_moment-angle,
 author = {Ioannis Gkeneralis},
 title = {Equivariant rigidity of complex and quaternionic moment--angle manifolds},
 year = {2026},
 howpublished = {Preprint, {arXiv}:2604.18455 [math.{AT}] (2026)},
 url = {https://arxiv.org/abs/2604.18455},
 arXiv = {arXiv:2604.18455}
}

@article{capovilla_li_loh,
 author = {Capovilla, Pietro and Li, Kevin and L{\"o}h, Clara},
 title = {Combination of open covers with {{\(\pi_1\)}}-constraints},
 fjournal = {Bulletin of the London Mathematical Society},
 journal = {Bull. Lond. Math. Soc.},
 issn = {0024-6093},
 volume = {57},
 number = {12},
 pages = {3886--3901},
 year = {2025},
 language = {English},
 doi = {10.1112/blms.70231},
 zbMATH = {8140786},
 Zbl = {1584.55004}
}

@article{JanuszkiewiczSwiatkowski,
 author = {Januszkiewicz, Tadeusz and {\'S}wi{\k{a}}tkowski, Jacek},
 title = {Commensurability of graph products},
 fjournal = {Algebraic \& Geometric Topology},
 journal = {Algebr. Geom. Topol.},
 issn = {1472-2747},
 volume = {1},
 pages = {587--603},
 year = {2001},
 language = {English},
 doi = {10.2140/agt.2001.1.587},
 url = {https://eudml.org/doc/122161},
 zbMATH = {1667106},
 Zbl = {0998.20029}
}

@article{GandiniRueping,
 author = {Gandini, Giovanni and R{\"u}ping, Henrik},
 title = {The {Farrell}-{Jones} conjecture for graph products},
 fjournal = {Algebraic \& Geometric Topology},
 journal = {Algebr. Geom. Topol.},
 issn = {1472-2747},
 volume = {13},
 number = {6},
 pages = {3651--3660},
 year = {2013},
 language = {English},
 doi = {10.2140/agt.2013.13.3651},
 zbMATH = {6216312},
 Zbl = {1291.18019}
}

@misc{gkeneralis-real,
 author = {Ioannis Gkeneralis},
 title = {Remarks on {Topological} {Rigidity} of {Real} {Moment}-{Angle} {Manifolds}},
 year = {2026},
 howpublished = {Preprint, {arXiv}:2604.15462 [math.{GT}] (2026)},
 url = {https://arxiv.org/abs/2604.15462},
 arXiv = {arXiv:2604.15462}
}

@article{KasprowskiLiLueck,
 author = {Kasprowski, Daniel and Li, Kevin and L{\"u}ck, Wolfgang},
 title = {{{\(K\)}}- and {{\(L\)}}-theory of graph products of groups},
 fjournal = {Groups, Geometry, and Dynamics},
 journal = {Groups Geom. Dyn.},
 issn = {1661-7207},
 volume = {15},
 number = {1},
 pages = {269--311},
 year = {2021},
 language = {English},
 doi = {10.4171/GGD/598},
 zbMATH = {7362829},
 Zbl = {1466.18008}
}

@article{Li2023,
 author = {Li, Kevin},
 title = {Amenable covers of right-angled {Artin} groups},
 fjournal = {Bulletin of the London Mathematical Society},
 journal = {Bull. Lond. Math. Soc.},
 issn = {0024-6093},
 volume = {55},
 number = {2},
 pages = {978--989},
 year = {2023},
 language = {English},
 doi = {10.1112/blms.12771},
 zbMATH = {7738032},
 Zbl = {1522.20140}
}

@article{AmelotteBriggs, title={Product decompositions of moment-angle manifolds and B-rigidity}, volume={66}, DOI={10.4153/S0008439523000383}, number={4}, journal={Canadian Mathematical Bulletin}, author={Amelotte, Steven and Briggs, Benjamin}, year={2023}, pages={1313–1325}}

@article{Wegner,
    author = {Wegner, Christian},
    title = {The Farrell–Jones conjecture for virtually solvable groups},
    journal = {Journal of Topology},
    volume = {8},
    number = {4},
    pages = {975-1016},
    year = {2015},
    month = {12},
    issn = {1753-8416},
    doi = {10.1112/jtopol/jtv026},
    url = {https://doi.org/10.1112/jtopol/jtv026},
    eprint = {https://academic.oup.com/jtopol/article-pdf/8/4/975/2760887/jtv026.pdf},
}

@article{BartelsLuck,
 author = {Bartels, Arthur and L{\"u}ck, Wolfgang},
 title = {The {Borel} conjecture for hyperbolic and {CAT}(0)-groups},
 fjournal = {Annals of Mathematics. Second Series},
 journal = {Ann. Math. (2)},
 issn = {0003-486X},
 volume = {175},
 number = {2},
 pages = {631--689},
 year = {2012},
 language = {English},
 doi = {10.4007/annals.2012.175.2.5},
 zbMATH = {6025000},
 Zbl = {1256.57021}
}

@article{Knopf,
 author = {Knopf, Svenja},
 title = {Acylindrical actions on trees and the {Farrell}-{Jones} conjecture},
 fjournal = {Groups, Geometry, and Dynamics},
 journal = {Groups Geom. Dyn.},
 issn = {1661-7207},
 volume = {13},
 number = {2},
 pages = {633--676},
 year = {2019},
 language = {English},
 doi = {10.4171/GGD/499},
 url = {nbn-resolving.org/urn:nbn:de:hbz:6-85209633373},
 zbMATH = {7096493},
 Zbl = {1471.18013}
}

@unpublished{WuElementaryAmenableFJC,
  author  = {Wu, Xiaolei},
  title   = {Elementary Amenable Groups and the Farrell--Jones Conjecture},
  year    = {2015},
  note    = {Preprint},
  url     = {https://wuxiaolei.org/notes/},
  urldate = {2026-10-16},
}

@article{Delzant,
 author = {Delzant, Thomas},
 title = {On the acylindrical accessibility of finitely presented groups},
 fjournal = {Annales de l'Institut Fourier},
 journal = {Ann. Inst. Fourier},
 issn = {0373-0956},
 volume = {49},
 number = {4},
 pages = {1215--1224},
 year = {1999},
 language = {French},
 doi = {10.5802/aif.1714},
 url = {https://eudml.org/doc/75378},
 zbMATH = {1323234},
 Zbl = {0999.20017}
}

@article{sela,
 author = {Sela, Z.},
 title = {Acylindrical accessibility for groups},
 fjournal = {Inventiones Mathematicae},
 journal = {Invent. Math.},
 issn = {0020-9910},
 volume = {129},
 number = {3},
 pages = {527--565},
 year = {1997},
 language = {English},
 doi = {10.1007/s002220050172},
 zbMATH = {1054509},
 Zbl = {0887.20017}
}

@article{Bowditch08,
 author = {Bowditch, Brian H.},
 title = {Tight geodesics in the curve complex},
 fjournal = {Inventiones Mathematicae},
 journal = {Invent. Math.},
 issn = {0020-9910},
 volume = {171},
 number = {2},
 pages = {281--300},
 year = {2008},
 language = {English},
 doi = {10.1007/s00222-007-0081-y},
 zbMATH = {5248088},
 Zbl = {1185.57011}
}

@article{gromov1982volume,
  title={Volume and bounded cohomology},
  author={Gromov, Michael},
  journal={Publications Math{\'e}matiques de l'IH{\'E}S},
  volume={56},
  pages={5--99},
  year={1982}
}

\end{document}